\documentclass[12pt]{article}

\title{Parametrized $\diamondsuit$ principles}
\author{Justin Tatch Moore, Michael Hru\v{s}\'{a}k, Mirna D\v{z}amonja
\footnote{
The first and third authors received support from EPSRC grant GR/M71121 for the
research of this paper.
The research of the second author was supported in part by the Netherlands 
Organization for Scientific Research (NWO) -- Grant 613.007.039, and in part 
by the Grant Agency of the Czech Republic -- Grant GA\v{C}R 201/00/1466.}
}
\usepackage{amsmath}
\usepackage{amsthm}
\usepackage{mathrsfs}
\usepackage{amssymb}
\usepackage{enumerate}

\newtheorem{thm}{Theorem}[section]
\newtheorem{lem}[thm]{Lemma}

\newtheorem{prop}[thm]{Proposition}

\newtheorem{cor}[thm]{Corollary}
\newtheorem*{thm*}{Theorem}
\newtheorem*{cor*}{Corollary}
\newtheorem*{question*}{Question}

\theoremstyle{definition}
\newtheorem{defn}[thm]{Definition}

\newtheorem*{notn}{Notation}

\newtheorem{question}[thm]{Question}

\theoremstyle{remark}
\newtheorem{remark}[thm]{Remark}

\newcommand\afrak{\mathfrak{a}}
\newcommand\bfrak{\mathfrak{b}}
\newcommand\cfrak{\mathfrak{c}}
\newcommand\dfrak{\mathfrak{d}}
\newcommand\ifrak{\mathfrak{i}}

\newcommand\rfrak{\mathfrak{r}}
\newcommand\sfrak{\mathfrak{s}}
\newcommand\tfrak{\mathfrak{t}}
\newcommand\ufrak{\mathfrak{u}}

\newcommand\Acal{\mathcal{A}}
\newcommand\Bcal{\mathcal{B}}
\newcommand\Ccal{\mathcal{C}}
\newcommand\Dcal{\mathcal{D}}
\newcommand\Ecal{\mathcal{E}}

\newcommand\Mcal{\mathcal{M}}
\newcommand\Ncal{\mathcal{N}}
\newcommand\Pcal{\mathcal{P}}
\newcommand\Qcal{\mathcal{Q}}
\newcommand\Rcal{\mathcal{R}}

\newcommand\Tcal{\mathcal{T}}
\newcommand\Ucal{\mathcal{U}}

\newcommand\Qbb{\mathbb{Q}}
\newcommand\Rbb{\mathbb{R}}
\newcommand\Sbb{\mathbb{S}}

\newcommand\add{\operatorname{add}}
\newcommand\cf{\operatorname{cf}}
\newcommand\cof{\operatorname{cof}}
\newcommand\non{\operatorname{non}}
\newcommand\cov{\operatorname{cov}}

\newcommand{\dom}{\operatorname{dom}}

\newcommand{\one}{\mathbf{1}}

\newcommand{\Th}{{}^{\textrm{th}}}

\newcommand{\Vbf}{\mathbf{V}}

\newcommand\axiom{\textrm}

\newcommand\MA{\axiom{MA}}
\newcommand\CH{\axiom{CH}}

\newcommand\ZFC{\axiom{ZFC}}

\newcommand\CON{\textrm{CON}}

\newcommand{\<}{\langle}
\renewcommand{\>}{\rangle}

\newcommand\N{\omega}
\newcommand\NN{\N^{\N}}

\newcommand\mand{\textrm{ and }}

\renewcommand\diamond{\diamondsuit}
\newcommand\club{\clubsuit}

\begin{document}

\maketitle

\begin{abstract}
We will present a collection of guessing principles which
have a similar relationship to $\diamond$ as cardinal
invariants of the continuum have to $\CH$.
The purpose is to provide a means for systematically analyzing
$\diamond$ and its consequences.
It also provides for a unified approach for understanding the status of a
number of consequences of $\CH$ and $\diamond$
in models such as those of Laver, Miller, and Sacks.
\end{abstract}

\section{Introduction}

Very early on in the course of modern set theory, Jensen
isolated the following combinatorial principle known as $\diamond$:
\begin{description}

\item[$\diamond$] There is a sequence $A_\alpha$ $(\alpha < \omega_1)$
such that for all $\alpha < \omega_1$, $A_\alpha \subseteq \alpha$ and
if $X$ is a subset of $\omega_1$ then set
$$\{\alpha < \omega_1:X \cap \alpha = A_\alpha\}$$
is stationary.

\end{description}
Jensen used this principle to construct a Suslin tree \cite{Souslin_Jensen}
and later many other constructions were carried out using
$\diamond$ as an assumption --- see \cite{more_settheory_top}.

The purpose of this paper is to provide a broad framework for analyzing
the consequences of Jensen's $\diamond$ principle.
Our intent is to present an array of $\diamond$-principles which have
the same relation to $\diamond$ as the cardinal invariants of the
continuum
(see e.g. \cite{set_theory_reals} or \cite{integers_topology}) have to
the Continuum Hypothesis.
We will approach an analysis of $\diamond$ in much the same
way as Blass approaches cardinal invariant inequalities
in \cite{reductions_characteristics}.

Our immediate motivation in this consideration stems from the isolation of
the principle $\diamond_{\dfrak}$ in \cite{thesis:Hrusak} (see also
\cite{another_diamond}).
\begin{description}

\item[$\diamond_{\dfrak}$]
There is a sequence $g_\alpha:\alpha \to \omega$ indexed by $\omega_1$
such that for every $f:\omega_1 \to \omega$ there is an $\alpha\geq
\omega$
with $f \restriction \alpha <^* g_\alpha$.

\end{description}
This principle implies $\dfrak = \omega_1$ (for much the same reason as
$\diamond$ implies $\cfrak = \omega_1$), follows
from $\dfrak = \omega_1 + \clubsuit$, and holds in most of the
standard generic extensions in which $\dfrak = \omega_1$ holds
\cite{another_diamond}.\footnote{
Another $\diamond$-like principle in this spirit is
the statement $\diamond(\omega_1^{<\omega})$ presented in
Section 6.2 of \cite{Pmax} (Definition 6.37).
To draw an analogy, this principle might also be called
$\diamond_{\non(\Mcal)}$ in the language of \cite{another_diamond}
(see Theorem 6.49 of \cite{Pmax}).
Also, Shelah has considered some specific cases of
$\Phi(A,B,E)$ defined below in the appendix of \cite{proper_forcing}.}
The main interest in it comes from the following fact relating to
the question of whether $\dfrak = \omega_1$ implies $\afrak = \omega_1$.
\begin{thm} \cite{another_diamond}\label{diamond_d/a}
$\diamond_{\dfrak}$ implies $\afrak = \omega_1$.
\end{thm}

It was initially unclear whether other cardinal invariants
of the continuum, such as $\bfrak$ and $\sfrak$, have similar
$\diamond$-like principles corresponding to them.
The cardinal $\sfrak$ was of particular interest in this context,
since it seemed that the construction of the Ostaszewski space
from $\omega_1$ random reals in \cite{thesis} should
be a consequence of a principle similar to $\diamond_{\sfrak}$
(whatever that might be).

It turned out that the correct language to use for formulating
$\diamond$-principles like those mentioned above was developed
by Devlin and Shelah in \cite{weak_diamond}.
They considered the following statement
\begin{description}

\item[$\Phi$]
For every $F:2^{<\omega_1} \to 2$ there is a $g:\omega_1 \to 2$
such that for every $f:\omega_1 \to 2$ the set
$\{\alpha \in \omega_1:F(f \restriction \alpha) = g(\alpha)\}$
is stationary.

\end{description}
which they showed to be equivalent to $2^{\aleph_0} < 2^{\aleph_1}$.
The framework of the weak diamond principle $\Phi$ of \cite{weak_diamond}
allows for the definition of two classes of
$\diamond$-principles, $\Phi(A,B,E)$ and $\diamond(A,B,E)$,
each taking a cardinal invariant $(A,B,E)$ as a parameter.

Like $\diamond_{\dfrak}$, these principles all imply that
the corresponding cardinal invariant is $\omega_1$.
They all follow from $\diamond$, with $\diamond(\cfrak)$ and
$\Phi(\cfrak)$ both being
equivalent to $\diamond$.
They also each have a ``guessing'' component which allows them to
carry out constructions for which one historically has used $\diamond$.
Moreover, many of the classical $\diamond$ constructions seem
to fit very naturally into this scheme.
For instance the standard construction of a Suslin tree from
$\diamond$ really requires only $\diamond(\non(\Mcal))$.
Also like $\diamond_{\dfrak}$, the principles
$\diamond(A,B,E)$ hold in many of the natural
models in which their corresponding cardinal invariant $\<A,B,E\>$
is $\omega_1$.
For instance $\diamond(\bfrak)$ holds in Miller's model and
$\diamond(\cof(\Ncal))$ holds in both the iterated and
the ``side-by-side'' Sacks models.

The paper is organized as follows.
Section \ref{Phi:sec}
introduces an abstract form of a cardinal invariant of the continuum
and formulates the principles $\Phi(A,B,E)$ which serve as
a first approximation to $\diamond(A,B,E)$.
Section \ref{Suslin_construction:sec}
presents the Suslin tree construction from
$\diamond$ in the language of $\Phi(\non(\Mcal))$.
Section \ref{diamond:sec} introduces a refinement of $\Phi(A,B,E)$
called
$\diamond(A,B,E)$ and gives some explanation for our choice of
it over $\Phi(A,B,E)$.
Section \ref{constructions:sec}
presents some more constructions which use $\diamond(A,B,E)$.
Section \ref{canonical_models:sec} shows that $\diamond(A,B,E)$ holds
in many of the models of $\<A,B,E\> = \omega_1$.
Section \ref{structured_invariants:sec}
studies the role of $\diamond(A,B,E)$ in
studying cardinal invariants other than those fitting into our
framework.
Section \ref{nonstandard_invariants:sec}
presents a proof that $\diamond(A,B,E)$ is not
a consequence of $\CH$ for any of the classical invariants $(A,B,E)$.

Our notation is, for the most part, standard (see \cite{set_theory}).
We will use $A^B$ to denote the collection of all functions from
$B$ to $A$.
$2^{<\omega_1}$ will be used to denote the tree of all functions
from a countable ordinal into $2$ ordered by extension.
If $t$ is a function defined on an ordinal, then we will
use $|t|$ to denote the domain of $t$.
Otherwise $|A|$ will be used to denote the cardinality of a set $A$.
The meaning of $|\cdot|$ should always be clear from the context.
If $B$ is a Borel subset of a Polish space, we will often identify it
with its code and use this code to define $B$ in forcing extensions.
We will use $\check B$ to represent the name for this set in the forcing
extension.
 
Many of the constructions in this paper will
require choosing a sequence $e_\delta:\omega \leftrightarrow \delta$
of bijections for each $\delta \in \omega_1$ or an increasing
sequence $\delta_n$ $(n \in \omega)$ which is cofinal in $\delta$ for
limit $\delta$.
To avoid repetition, we will fix a sequence of bijections
$e_\delta$ $(\delta \in \omega_1)$ and
cofinal sequences $\delta_n$ $(n \in \omega)$ for limit $\delta$
once and for all.
If there is a need to refer to, e.g., a special cofinal
sequence in $\delta$ we will use $\bar \delta_n$ $(n \in \omega)$ for the
sequence instead.

\section{Abstract cardinal invariants and $\Phi$}
\label{Phi:sec}

The following structure allows for a compact definition
of many of common cardinal invariants of the continuum.

\begin{defn} \cite{Vojtas}
An \emph{invariant} is a triple
$(A,B,E)$ such that
\begin{enumerate}

\item $A$ and $B$ are sets of cardinality at most $|\Rbb|$,

\item $E \subseteq A \times B$,

\item for every $a\in A$ there is a $b\in B$ such that
$(a,b)\in E$,

\item for every $b\in B$ there is an $a \in A$ such that
$(a,b)\not\in E$.

\end{enumerate}
\end{defn}

Usually we will write $a E b$ instead of $(a,b)\in E$.

\begin{defn}
If $(A,B,E)$ is an invariant, then its \emph{evaluation}
$\<A,B,E\>$ is given by
$$\<A,B,E\> = \min\{|X|:X \subseteq B \mand (\forall a \in A)( \exists b
\in X)
(a E b)\}$$
\end{defn}

If $A = B$ then we will write $(A,E)$ and $\<A,E\>$ instead of $(A,B,E)$
and
$\<A,B,E\>$ respectively.
Two typical examples of invariants are
$(\Ncal,\subseteq)$ and
$(\Mcal,\Rbb,\not \ni)$.
The evaluations $\<\Ncal,\subseteq\>$ and $\<\Mcal,\Rbb,\not \ni\>$
are clearly just $\cof(\Ncal)$ and $\non(\Mcal)$. Even though, strictly
speaking,
$\Mcal$ and $\Ncal$ are ideals of cardinality $2^\cfrak$ they both have
a basis consisting
of Borel sets, hence of cardinality $\cfrak$.  
If an invariant $(A,B,E)$ already has a common representation,
we will use such a representation instead of $(A,B,E)$.
Moreover, we will abuse notation and use these representations to
abbreviate both the invariant and its evaluation.
What we mean should always be clear from the context.

\begin{defn}
Let $(A,B,E)$ be an invariant.
$\Phi(A,B,E)$ is the following statement:
\begin{description}

\item[$\Phi(A,B,E)$]
For every $F:2^{<\omega_1} \to A$ there is a $g:\omega_1 \to B$
such that for every $f:\omega_1 \to 2$ the set
$\{\alpha \in \omega_1:F(f \restriction \alpha) E\, g(\alpha)\}$
is stationary.

\end{description}
The witness $g$ for a given $F$ in this statement will be called
a \emph{$\diamond(A,B,E)$-sequence for $F$}.
If $F(f \restriction \delta) E g(\delta)$ then we will say that 
$g$ guesses $f$ (via $F$) at $\delta$.
\end{defn}

\begin{prop}
\label{diamond_implies}
$\diamond$ implies $\Phi(A,B,E)$ for any invariant $(A,B,E)$.
\end{prop}

\begin{proof}
Let $A_\alpha$ $(\alpha \in \omega_1)$ be a diamond sequence which
guesses elements of $2^{\omega_1}$
($A_\alpha$ is in $2^\alpha$).
Set $g(\alpha)$ to be any $b\in B$ such that
$F(A_\alpha) E b$.
Then $g$ is a $\diamond(A,B,E)$-sequence for $F$ since for all
$f:\omega_1 \to 2$
$$\{\delta \in \omega_1:f \restriction \delta = A_\delta\}
\subseteq \{\delta \in \omega_1:F(f \restriction \delta) E g(\delta)\}.$$
\end{proof}

\begin{prop}
$\Phi(A,B,E)$ implies $\<A,B,E\>$ is at most $\omega_1$.
\end{prop}

\begin{proof}
Let $F:2^\omega \to A$ be a surjection and extend
$F$ to $2^{<\omega_1}$ by setting
$F(t) = F(t \restriction \omega)$ if $t$ has an infinite domain
and defining $F(t)$ arbitrarily otherwise.
Let $g$ be a $\diamond(A,B,E)$-sequence for $F$.
It is easy to see that the range of $g$ witnesses
$\<A,B,E\> \leq \omega_1$.
\end{proof}

Notice the resemblance of this proof to the standard proof that
$\diamond$ implies $\CH$.
In fact, if we view $\cfrak$ as the invariant $(\Rbb,=)$ then
we have the following fact.

\begin{prop}
\label{Phi_c}
$\Phi(\cfrak)$ is equivalent to $\diamond$.
\end{prop}

\begin{proof}
By Proposition \ref{diamond_implies}
we need only to show that
$\Phi(\cfrak)$ implies $\diamond$.
For each infinite
$\alpha \in \omega_1$, fix a bijection $H_\alpha:2^\alpha \to \Rbb$.
Set $F(t) = H_\alpha(t)$ where $\alpha = |t|$.
Let $g$ be the $\diamond(\cfrak)$-sequence for this $F$.
Set $A_\alpha = H^{-1}_\alpha(g(\alpha))$.
It is easy to see that $A_\alpha$ $(\alpha \in \omega_1)$ is a
$\diamond$-sequence.
\end{proof}

\begin{prop}
$\Phi(\Rbb^\omega,\sqsubseteq)$ is equivalent to $\diamond$.
Here $f \sqsubseteq g$ iff the range of $f$ is contained in the
range of $g$.
\end{prop}

\begin{proof}
Combine the previous proof with the Kunen's result stating
that $\diamond$ is equivalent to $\diamond^-$ (see \cite{set_theory}).
\end{proof}

A natural question which arises is:
``When do relations between
invariants translate into implications between the corresponding
$\diamond$-principles?''
This is largely answered by the next proposition.

\begin{notn} (\emph{Tukey ordering} \cite{Vojtas})
If $(A_1,B_1,E_1)$ and $(A_2,B_2,E_2)$ are invariants then
$$(A_1,B_1,E_1) \leq_T (A_2,B_2,E_2)$$ 
when there are maps
$\phi:A_1 \to A_2$ and $\psi:B_2 \to B_1$ such that
$(\phi(a),b)\in E_2$ implies $(a,\psi(b))\in E_1$.
\end{notn}

As one would expect, the Tukey ordering on invariants gives
the corresponding implications for $\Phi$ principles.

\begin{prop}
\label{Phi_Tukey}
If $(A_1,B_1,E_1) \leq_T (A_2,B_2,E_2)$ then $\Phi(A_1,B_1,E_1)$
is a consequence of $\Phi(A_2,B_2,E_2)$.
\end{prop}

One should exercise caution, however, when trying to turn inequalities
between evaluations of cardinal invariants into implications between
$\diamond$-principles.
For instance, $(\omega^\omega,<^*)$ are $(\omega^\omega,\leq)$ have
the same evaluation but seem to give rise to different
$\diamond$-principles.
We will use $\dfrak$ to denote $(\omega^\omega,<^*)$.

The smallest invariant in the Tukey order is $(\Rbb,\ne)$.
It is known that $\Phi(2,\ne)$ is equivalent to $\Phi(\Rbb,\ne)$ ---
this was noted by Abraham and can be extracted from
\cite{weak_diamond}.
The proof is given for completeness.

\begin{thm}
$\Phi(2,\ne)$ is equivalent to $\Phi(\Rbb,\ne)$.
\end{thm}

\begin{proof}
Since $(\Rbb,\ne)$ is below $(2,\ne)$ in the Tukey
order, it suffices to show that $\Phi(\Rbb,\ne)$ implies
$\Phi(2,\ne)$.
To this end, suppose that $F:2^{<\omega_1} \to 2$ witnesses that
$\Phi(2,\ne)$ fails.
Define a function $F^*$ whose ranges is contained in $2^\omega$ and
whose domain consists of functions of the form
$t:\delta \times \omega \to 2$ so that
$F^*(t)(i) = F(t(\cdot,i))$.
Now let $g:\omega_1 \to 2^\omega$ be given.
To see that $g$ is not a $\diamond(\Rbb,\ne)$-sequence for $F^*$,
pick closed unbounded sets $C_n \subseteq \omega_1$ and
functions $f_n:\omega_1 \to 2$ such that
$F(f_n \restriction \delta) = g(\delta)(n)$
for every $n$ and $\delta$ in $C_n$.
Now define $f:\omega_1 \times \omega \to 2^\omega$ by putting
$f(\delta,n) = f_n(\delta)$.
Then $F^*(f \restriction \delta \times \omega) = g(\delta)$
whenever $\delta$ is in $\bigcap_{n \in \omega} C_n$.
\end{proof}

\section{The Suslin tree construction}

\label{Suslin_construction:sec}

In order to get a feel for how the statements $\Phi(A,B,E)$ are used,
we will begin by revisiting an old construction and translating it
into the language which we have developed.

\begin{thm} \label{Suslin_construction}
$\Phi(\non(\Mcal))$ implies that there is a Suslin tree.
\end{thm}

\begin{proof}
By some suitable coding, $F$ will take triples $(\alpha,\prec,A)$
as its argument
where $\alpha \in \omega_1$, $\prec \subseteq \alpha^2$,
and $A \subseteq \alpha$.
$F$ will be defined to be the empty set unless
\begin{enumerate}

\item $\alpha$ is a limit of limit ordinals,

\item
$\prec$ is a tree order on $\alpha \in \omega_1$ of limit height,

\item
if $\gamma < \alpha$ then for every $\delta$ less than the height of
$(\alpha,\prec)$ there
is a $\bar \gamma < \alpha$ with $\gamma \prec \bar \gamma$ and
the height of $\bar \gamma$ greater than $\delta$,

\item
every element of $\alpha$ has exactly $\omega$ immediate successors
in $\prec$,

\item
if $\xi < \alpha$,
then $[\xi \cdot \omega,\xi \cdot \omega + \omega)$ is exactly the
collection
of elements of $(\alpha,\prec)$ of height $\xi$, and

\item
$A$ is a maximal antichain in $(\alpha,\prec)$.

\end{enumerate}
For such a triple $(\alpha,\prec,A)$ let
$\alpha_n$ $(n \in \omega)$
be an increasing sequence cofinal in $\alpha$ and such that each
$\alpha_n$
is a limit ordinal.
Let $[\alpha,\prec]$ denote the collection of all cofinal branches
through $(\alpha,\prec)$.
Define $$\phi_\prec:[\alpha,\prec] \to \omega^\omega$$ by
setting $\phi_\prec(b)(n)$ to be the unique $k$ such that
$\alpha_n + k$ is in $b$.
Notice that if $A \subseteq \alpha$ is a maximal antichain
then
$$N(\alpha,\prec,A) =
\{\phi_\prec(b): b \in [\alpha,\prec] \mand A \cap b =\emptyset\}$$ is
closed and nowhere dense in $\NN$.
Also, observe that $\phi_\prec$ is a surjection.
Let $F(\alpha,\prec,A)$ be the collection of all finite changes of elements
of $N(\alpha,\prec,A)$.

Now suppose that $g:\omega_1 \to \NN$ is a $\diamond(\non(\Mcal))$-sequence for
$F$.
Construct a tree order $\prec$ on $\omega_1$ by recursion.
Define $(\omega^2,\prec)$
so that it is isomorphic to $\omega^{<\omega}$ ordered by end extension.
Now suppose that $(\alpha,\prec)$ is defined and has limit height.
Extend the order to $(\alpha + \omega,\prec)$ in such a way that
a cofinal branch $b$ in $[\alpha,\prec]$ has an upper bound
in $(\alpha+\omega,\prec)$ iff
$\phi_\prec(b)$ is eventually equal to $g(\alpha)$.
Now extend $\prec$ to $\alpha + \omega^2$ in such a way that
conditions 1-5 above are satisfied.

To see that $(\omega_1,\prec)$ is a Suslin tree, suppose that
$A \subseteq \omega_1$ is a maximal antichain in $(\omega_1,\prec)$.
By the same crucial lemma
as in the standard $\diamondsuit$ construction
(see Lemma 7.6 in Chapter II of \cite{set_theory})
the set of $\alpha<\omega_1$ such that $A\cap\alpha$ is a maximal
antichain in $(\alpha,\prec)$
contains a closed unbounded set $C$. If $g$ guesses $(\omega_1,\prec,A)$
at $\alpha \in C$
and $F(\alpha,\prec \restriction \alpha,A \cap \alpha)$
is nonempty, then $A \cap \alpha$ is a maximal antichain in
$(\alpha + \omega,\prec)$.
It is now easily verified using properties 1--4 above that,
since $A \subseteq \alpha$ is a maximal antichain in $(\alpha +
\omega,\prec)$,
it is maximal in $(\omega_1,\prec)$ as well.
\end{proof}

A reader familiar with the classical construction of a Suslin tree
from $\diamond$ (see, e.g., Section II.7 of \cite{set_theory})
should have no trouble in seeing that this is indeed the
same construction with the assumption reduced to the minimum required to
carry out the argument.
In Section \ref{diamond:sec} we shall comment that the principle
$\diamond(\non(\Mcal))$ implied by $\Phi(\non(\Mcal))$ suffices for this
construction and in Section
\ref{canonical_models:sec} we will
see that $\diamond(\non(\Mcal))$ is in fact
much weaker than $\diamond$.

\section{$\diamond(A,B,E)$ --- a definable form of $\Phi(A,B,E)$}
\label{diamond:sec}

The purpose of this section is to demonstrate that $\Phi(A,B,E)$ is,
in general, too strong to hold in typical generic extensions
in which $\<A,B,E\> = \omega_1$.
We will, however, recover from it a principle $\diamond(A,B,E)$ which is
still strong enough for most of the combinatorial applications of
$\Phi(A,B,E)$ and which is of a more appropriate strength.

\begin{prop}
\label{Phi_gives_WD}
$\Phi(A,B,E)$ implies $2^{\omega} < 2^{\omega_1}$.
\end{prop}

\begin{proof}
Suppose that $2^{\omega} = 2^{\omega_1}$. 
Let $H:2^\omega \to B^{\omega_1}$ be a surjection.
Define $F(t)$ to be any $a$ in $A$ such that
$(a,H(t \restriction \omega)(|t|))$ is not in $E$.
Now if $g:\omega_1 \to B$ is given, pick an $f:\omega_1 \to 2$
such that $H(f \restriction \omega) = g$.
It is easily checked that $g$ does not guess $f$ at any
$\delta \geq \omega$.
\end{proof}

A closer look at the uses of $\Phi(A,B,E)$ presented in Sections
\ref{Phi:sec} and \ref{Suslin_construction:sec}
reveals that in all cases the maps $F$ which were used in the proofs
could be chosen to be nicely definable.
This, generally speaking, is atypical of the map
$F$ in the proof of Proposition \ref{Phi_gives_WD}.
Before we discuss the principles $\diamond(A,B,E)$, we first
need to define the notion of a Borel invariant.

\begin{defn} \cite{reductions_characteristics}
An invariant $(A,B,E)$ is \emph{Borel} if $A$, $B$ and $E$ are Borel
subsets of some Polish space.
\end{defn}

With slight technical changes, all of the ``standard''
invariants $(A,B,E)$
can be represented as Borel invariants.
The invariants for which this is non trivial are those in
Cicho\'n's diagram.
First note that the $G_\delta$ null and $F_\sigma$ meager sets generate
$\Ncal$ and $\Mcal$ respectively.
Furthermore, $\subseteq$ is a Borel relation on a cofinal subset
of the $G_\delta$ null and $F_\sigma$ meager sets. 
The details for the category invariants are handled in
Section 3 of \cite{reductions_characteristics}.
For null sets, one can use the fact that any null set
is contained in the union of two \emph{small} sets and that the containment
relation on such unions is Borel
(see Section 2.5 of \cite{set_theory_reals} for a discussion of small set
and their relation to null sets).

\begin{defn}
Suppose that $A$ is a Borel subset of some Polish space.
A map $F:2^{<\omega_1} \to A$ is \emph{Borel} if for every
$\delta$ the restriction of $F$ to $2^\delta$ is a Borel map.
\end{defn}

As we will see throughout this paper, the maps $F$ which we are
actually interested in considering in the context of $\Phi(A,B,E)$
all can be made to satisfy this requirement.
This motivates the following definition.

\begin{defn}
Let $(A,B,E)$ be a Borel invariant.
$\diamond(A,B,E)$ is the following statement:
\begin{description}

\item[$\diamond(A,B,E)$]
For every Borel map $F:2^{<\omega_1} \to A$ there is a $g:\omega_1 \to B$
such that for every $f:\omega_1 \to 2$ the set
$\{\alpha \in \omega_1:F(f \restriction \alpha) E g(\alpha)\}$
is stationary.

\end{description}
\end{defn}

Aside from the fact that $\diamond(A,B,E)$ often suffices for
applications
of $\Phi(A,B,E)$, it is also the case that, unlike $\Phi(A,B,E)$,
$\diamond(A,B,E)$ is often forced in the standard models where
$\<A,B,E\> = \omega_1$ is forced.
This is the content of
Section \ref{canonical_models:sec}.
The key property of Borel maps which we will need in
Section \ref{canonical_models:sec} is that
if $M$ is a model of $\ZFC$ (usually an intermediate forcing
extension) which contains the codes for $A$ and $F \restriction
2^\delta$
and $t \in 2^{\delta}$
then $F(t)$ can be computed in $M$.
Often it will be convenient to define a map $F$ only on
a Borel subset of $2^\delta$ for each $\delta$.
In such a case $F$ will assume a fixed constant value elsewhere.

The reader is now encouraged to re-read Section \ref{Phi:sec} and
convince themselves that $\diamond(A,B,E)$ suffices in each case
in which $\Phi(A,B,E)$ was used as an assumption
for a particular Borel invariant $(A,B,E)$.
For instance we have the following theorems.

\begin{prop}
$\diamond(\cfrak)$ is equivalent to $\diamond$.
\end{prop}

\begin{prop}
\label{diamondR/diamond2}
$\diamond(\Rbb,\ne)$ is equivalent to $\diamond(2,\ne)$.
\end{prop}

\begin{thm}
$\diamond(\non(\Mcal))$ implies the existence of a Suslin tree.
\end{thm}

Another problematic aspect of the statements
$\Phi(A,B,E)$ is that under $\CH$, the Tukey types of many of the
standard invariants are reduced to $(\omega_1,<)$.
For instance, under $\dfrak = \omega_1$
the Tukey type of $(\omega^\omega,<^*)$ reduces to
$(\omega_1,<)$ and hence
$\Phi(\dfrak)$ is equivalent to $\Phi(\omega_1,<)$
(see \cite{thesis:Yiparaki}).
The Tukey maps in such situations,
however, are generally far from being definable.
The analog of Proposition \ref{Phi_Tukey} for $\diamond(A,B,E)$
avoids this.

\begin{defn}\label{Borel_Tukey}
(Borel Tukey ordering \cite{reductions_characteristics})
Given a pair of Borel invariants $(A_1,B_1,E_1)$,  $(A_2,B_2,E_2)$
We will say that $(A_1,B_1,E_1) \leq_T^B (A_2,B_2,E_2)$ if
$(A_1,B_1,E_1) \leq_T (A_2,B_2,E_2)$ and the connecting maps are
both Borel.
\end{defn}

\begin{prop}
\label{diamond_Tukey}
If $(A_1,B_1,E_1) \leq^B_T (A_2,B_2,E_2)$ then
$\diamond(A_2,B_2,E_2)$ implies $\diamond(A_1,B_1,E_1)$.
\end{prop}

It turns out that the Tukey connections between all the invariants
we will consider satisfy the above
requirement (see \cite{reductions_characteristics})
and hence implications such as
$\diamond(\add(\Mcal))$ implies $\diamond(\add(\Ncal))$ hold.

\section{Some more constructions}
\label{constructions:sec}

The purpose of this section is to present some more topological and
combinatorial
constructions.
The first construction is that of the Ostaszewski space of
\cite{Ostaszewski_space}.
Recall that an Ostaszewski space is a countably compact non-compact
perfectly normal space.
Usually this space is considered to have the additional property that
its closed sets are either countable or co-countable.
Originally this space was constructed using $\clubsuit + \CH$, an
equivalent of $\diamond$ \cite{Ostaszewski_space}.

Unlike the example of the Suslin tree, which does not seem to
yield any new models in which there are Suslin trees,
the hypothesis we use in the construction
makes it rather transparent that there are Ostaszewski spaces after
adding $\omega_1$ random reals.
The construction of an Ostaszewski space from a sequence of random reals
(see \cite{thesis} or \cite{random_top}) and a careful analysis of the
combinatorics involved
was one of the main motivations and inspirations for the formulation
of $\diamond(\sfrak)$
and consequently $\diamond(A,B,E)$ for arbitrary invariants $(A,B,E)$.
  
\begin{notn}
Let $(\omega)^\omega_\omega$ denote the collection of all partitions
of $\omega$ into infinitely many infinite pieces.
\end{notn}

Recall that if $A,B \subseteq \omega$ then $A$ is \emph{split by $B$} if
both $A \cap B$ and $A \setminus B$ are infinite.
The invariant which seems to be at the heart of Ostaszewski's
construction is
$$\sfrak^\omega = ([\omega]^\omega,(\omega)^\omega_\omega,
\textrm{is split by all pieces of}),$$
a close relative of
$$\sfrak = ([\omega]^\omega,\textrm{is split by}).$$
The invariant $\sfrak^\omega$ is connected to $\non(\Mcal)$ and
$\non(\Ncal)$ by the following proposition.

\begin{prop}
$\sfrak^\omega$ is below both $\non(\Mcal)$ and $\non(\Ncal)$ in
the Borel Tukey order.
\end{prop}

\begin{proof}
Let $\mu$ be the product measure on $\omega^\omega$ obtained by
setting $\mu(\{n\}) = 2^{-n-1}$.
Let $S$ be the collection of all $f$ in $\omega^\omega$ which take
the value $n$ infinitely often for each $n$.
It is easily verified that $S$ is both comeager and measure 1 and hence we
can view $\non(\Mcal) = (\Mcal,S,\not \ni)$ and
$\non(\Ncal) = (\Ncal,S,\not \ni)$.
Define $\phi:[\omega]^\omega \to \Mcal \cap \Ncal$ by letting
$\phi(A)$ be the collection of all $f$ in $S$ which take
all values infinitely often on $A$.
Define $\psi:S \to (\omega)^\omega_\omega$ by
$\psi(f) = \{f^{-1}(n):n \in \omega\}$.
It is easily verified that this pair of maps gives the desired
Borel Tukey connections.
\end{proof}

\begin{thm} \label{Ostaszewski_construction}
$\diamond(\sfrak^\omega)$ implies the existence of a perfectly normal
countably compact non-compact space (i.e. an Ostaszewski space).
\end{thm}

\begin{proof}
Again by suitable coding, we will take the domain of $F$ to
be the set of all triples $(\alpha,\Bcal,D)$ such that
$\alpha \in \omega_1$, $\Bcal = \<U_\gamma:\gamma < \alpha\>$ where
$U_\gamma \subseteq \gamma + 1$, $\gamma \in U_\gamma$,
and $D \subseteq \alpha$.
Given a pair $(\alpha,\Bcal)$ as above, let $\tau_{\Bcal}$ be the
topology on $\alpha$ generated by taking $\Bcal$ as a clopen subbase.
$F(\alpha,\Bcal,D)$ is defined to be $\omega$ unless
\begin{enumerate}

\item $\alpha$ is a limit ordinal,

\item $U_\gamma$ is compact in $\tau_{\Bcal}$ for all $\gamma<\alpha$,

\item \label{final_segment}
for every $\gamma < \alpha$, the closure of $[\gamma,\gamma+\omega)$
in $(\alpha,\tau_{\Bcal})$ is $[\gamma,\alpha)$, and 

\item $D$ does not have compact closure in $(\alpha,\tau_{\Bcal})$.

\end{enumerate}
Define $V_{\alpha,n}$ for $n$ in $\omega$ by setting
$V_{\alpha,0} = U_{e_{\alpha}(0)}$ and
$$V_{\alpha,n} = U_{e_\alpha(k)} \setminus \bigcup_{i < n} V_{\alpha,i}$$
where $k$ is minimal such that this set is nonempty and such that
$e_{\alpha}(n)$ is covered by $V_{\alpha,i}$ for some $i\leq n$.
Thus $\{V_{\alpha,n}:n \in \omega\}$ is a partition of
$(\alpha,\tau_{\Bcal})$ into compact open sets.
Set $F(\alpha,\Bcal,D)$ to be the collection of all $n$ such that
$D \cap V_{\alpha,n}$ is nonempty.
Notice that since $V_{\alpha,n}$ is compact for all
$n$ and $D$ does not have compact
closure, $F(\alpha,\Bcal,D)$ is infinite.

Now let $g:\omega_1 \to (\omega)^\omega_\omega$ be a
$\diamond(\sfrak^\omega)$-sequence
for $F$.
Define a locally compact
topology $(\omega_1,\tau_{\Bcal})$ by recursion.
Suppose that $\Bcal \restriction \alpha$ have been defined so far,
satisfying 1-3.
Notice that if $A \subseteq \omega$ and $\sigma,\sigma' \in
(\omega)^\omega_\omega$
are such that $A$ is split by all pieces of $\sigma$ and every element
of
$\sigma'$ contains some element of $\sigma$ then $A$ is split by
every element of $\sigma'$.
Thus by altering $g(\alpha)$, if necessary, we may assume that for
each $k$, the collection $\{V_{\alpha,n}:n \in g(\alpha)(k)\}$ has a
union
which is cofinal in $\alpha$.
Let $$U_{\alpha + k} = \bigcup_{n \in g(\alpha)(k)} V_{\alpha,n}.$$
Since $U_{\alpha+k}$ is cofinal in $\alpha$ for all $k$,
the closure of a co-bounded subset of $\alpha$ is co-bounded in
$\alpha + \omega$.

Clearly $(\omega_1,\tau_{\omega_1})$ is not compact since all
initial segments are open in $\tau_{\omega_1}$.
To finish the proof, it suffices to show that closed sets
are either compact or co-countable.
Now suppose that $D \subseteq \omega_1$ does not have compact closure.
Let $\delta \in \omega_1$ be such that
$F(\delta,\Bcal \restriction \delta,D \cap \delta)$ is defined and is
split by every element of $g(\delta)$.
Then $D \cap \delta$ must accumulate at $\delta + n$ for all $n$.
It follows from \ref{final_segment} that the closure of $D$ in
$\omega_1$ is co-bounded.
\end{proof}

As mentioned above, the construction can be carried out using
$\diamond(\non(\Ncal))$ which holds after adding
$\omega_1$ random reals.
Eisworth and Roitman have shown that the construction of
an Ostaszewski space can not be carried out under $\CH$ alone
and hence some form of a guessing principle is required
\cite{no_ostaszewski}.
While the space above is hereditarily separable,
the following question is open:

\begin{question}
Does $\diamond(\non(\Ncal))$ imply the existence of a non-metric compact
space
$X$ such that $X^2$ is hereditarily separable?
\end{question}

We will now pass to a purely combinatorial construction.
Recall that a sequence $A_\alpha:\alpha \to 2$ indexed by $\omega_1$ is
\emph{coherent} if for every $\alpha < \beta$
$A_\alpha =^* A_\beta \restriction \alpha$.
Such a sequence is \emph{trivial} if there is a $B:\omega_1 \to 2$ such
that $A_\alpha =^* B \restriction \alpha$ for all $\alpha$.
Non-trivial coherent sequences can be constructed without additional
set theoretic assumptions
\cite{properties_Stone_Cech}.
The conclusion of the following theorem is deduced from $\diamond$ in
\cite{more_settheory_top}, though unlike
Theorems \ref{Suslin_construction} and \ref{Ostaszewski_construction},
the argument presented here does not
mirror an existing argument.

\begin{thm} \label{big_sequence}
$\diamond(\bfrak)$ implies that there is a coherent sequence
$A_\alpha$ $(\alpha \in \omega_1)$ of binary maps
such that for every uncountable
set $X$, there is an $\alpha \in \omega_1$ with
$A_\alpha$ taking both its values infinitely often on $X \cap \alpha$.
\end{thm}

First we will need the following fact which seems to be
of independent interest.
Recall that a \emph{ladder system} is a sequence
$\<C_\delta:\delta \in \lim(\omega_1)\>$ such that
$C_\delta$ is a cofinal subset of
$\delta$ of order type $\omega$ for each limit ordinal $\delta \in \omega_1$.

\begin{thm} \label{weak_club}
$\diamond(\bfrak)$ implies that there is a ladder system
$C_\delta$ such that for every sequence
of uncountable sets $X_\gamma \subseteq \omega_1$ $(\gamma \in \omega_1)$
there
are stationarily many $\delta$ such that $X_\gamma \cap C_\delta$ is infinite
for all $\gamma < \delta$.
\end{thm}

\begin{proof} (of Theorem \ref{weak_club})
Let $\vec{X} = \<X_\gamma : \gamma < \delta\>$
be a given sequence of subsets of $\delta$
and set $F(\vec{X})$ to be the identity function unless $\delta$ is
a limit ordinal and $X_\gamma$ is unbounded in $\delta$
for all $\gamma < \delta$.
Set $$F(\vec{X})(n) =
\max_{i \leq n} \left\{\min
\{e^{-1}_\delta(\gamma):\gamma \in X_{e_{\delta}(i)} \setminus \delta_n\}
\right\}.$$

Now suppose that $g:\omega_1 \to \omega^\omega$ is a $\diamond(\bfrak)$-sequence
for $F$.
By making $g(\delta)$ larger if necessary, we may assume that
$$g(\delta)(n) > e^{-1}_{\delta}(\delta_n).$$
Set $$C_\delta = \bigcup_{n = 0}^\infty
\{\gamma < \delta:(e^{-1}_\delta(\gamma) \leq g(\delta)(n)) \mand
(\gamma \geq \delta_n)\}.$$
Clearly $C_\delta$ is an $\omega$-sequence which is cofinal in $\delta$
and
it is routine to check that it satisfies the conclusion of the theorem.
\end{proof}

\begin{defn}
A binary coherent sequence $\vec{A}$ \emph{almost contains}
a ladder system $\vec{C}$ if 
$A_\alpha$ is eventually 1 on $C_\delta$ whenever $\delta < \alpha$.
\end{defn}

Notice that if $\vec{C}$ satisfies the conclusion of Theorem
\ref{weak_club}
then for any binary sequence $A_\alpha$ which almost contains $\vec{C}$
and
any uncountable set $X$, there is an $\alpha$ such that
$A_\alpha$ takes the value 1 infinitely often on $X \cap \alpha$.
Since coherence implies that this occurs for all $\beta \geq \alpha$ as
well, it suffices to prove the following lemma.

\begin{lem}
$\diamond(\bfrak)$ implies that for every ladder system $\vec{C}$
there is
a coherent sequence $\vec{A}$ which almost contains $\vec{C}$ such that
for every uncountable set $X \subseteq \omega_1$ there is an $\alpha$
such that $A_\alpha$ takes the value 0 infinitely often on $X \cap
\alpha$.
\end{lem}

\begin{proof}
First recall the following notion of a minimal walk
(see \cite{cseq} or \cite{partitioning_ordinals}).
If $\alpha < \beta$ then $\beta(\alpha) = \min (C_\beta \setminus
\alpha)$.
Here $C_{\alpha+1} = \{\alpha\}$.
Define $\beta^i(\alpha)$ recursively by setting
$\beta^0(\alpha) = \beta$ and $\beta^{i+1}(\alpha) =
\beta^i(\alpha)(\alpha)$.
Let $$a_\beta(\alpha) = |C_{\beta^{k-1}(\alpha)} \cap \alpha|$$
where $k$ is the minimal such that $\beta^k(\alpha) = \alpha$.
That is, $a_\beta(\alpha)$ is the
weight of the last step in the walk from $\beta$ to $\alpha$.

Now let $X \subseteq \delta$ be given.
Define $F(X,\delta)$ to be the identity function unless $X$
is cofinal in $\delta$ in which case
set $$F(X,\delta)(n) = \min\{a_\delta(\gamma):\gamma \in \delta \cap X
\setminus \bar \delta_n\},$$
where $\bar \delta_n$ $(n \in \omega)$
is an increasing enumeration of $C_\delta$.

Let $g:\omega_1 \to \omega^\omega$ be a $\diamond(\bfrak)$-sequence for $F$.
By making functions in $g$ larger if necessary we may assume that
$g(\alpha)$ is monotonic for all $\alpha$.
Set $$b_\beta(\alpha) = 
\max_{i < k -1} g(\beta^i(\alpha))
(|C_{\beta^i(\alpha)} \cap \alpha|)$$ if
the maximum is over a nonempty set and 0 otherwise
(where, again, $k$ is minimal such that $\beta^k(\alpha) = \alpha$).
Define $A_\beta(\alpha)$ to be 0 if $a_\beta(\alpha) < b_\beta(\alpha)$
and 1 otherwise.

It is routine to show that $\vec{a}$ and $\vec{b}$
are both coherent and hence that $\vec{A}$ is coherent
(see section 1 of \cite{cseq} or \cite{partitioning_ordinals}).
It is equally routine to show that $b_\beta$ is eventually constant
on any ladder while $a_\beta$ is eventually 1-1 on each ladder and
hence $\vec{A}$ almost contains $\vec{C}$. 
To see that $\vec{A}$ satisfies the conclusion of the theorem,
let $X$ be an uncountable set.
Fix a $\delta$ such that $X \cap \delta$ is cofinal in $\delta$ and
$g(\delta)$ is not dominated by $F(X \cap \delta,\delta)$.
Let $\gamma < \delta$ be arbitrary.
It suffices to find an $\alpha$ in $X \cap \delta \setminus \gamma$
such that $A_\delta(\alpha) = 0$.
Let $n$ be such that $\bar \delta_n > \gamma$ and
$F(X \cap \delta,\delta)(n) < g(\delta)(n)$.
Let $\alpha \in X \cap \delta \setminus \bar \delta_n$ be such that
$a_\delta(\alpha) = F(X \cap \delta,\delta)(n)$.
Now $$b_\delta(\alpha) \geq g(\delta)(|C_\delta \cap \alpha|) \geq
g(\delta)(n) > a_\delta(\alpha)$$
which finishes the proof.
\end{proof}

\begin{question}
For which Borel invariants $(A,B,E)$ does
$\diamond(A,B,E)$ imply the existence of a c.c.c.
destructible $(\omega_1,\omega_1^*)$-gap in $[\omega]^\omega$?
\end{question}

\section{Canonical models for $\diamond(A,B,E)$}
\label{canonical_models:sec}

The purpose of this section is to show that, for the classical
invariants
$(A,B,E)$,
$\diamond(A,B,E)$ holds in many of the standard models for
$\<A,B,E\> = \omega_1$.

\begin{thm}
\label{Cohen_random_diamond}
Let $\Ccal_{\omega_1}$ and $\Rcal_{\omega_1}$ be the Cohen and measure
algebras
corresponding to the product space $2^{\omega_1}$ with its usual
topological
and measure theoretic structures.
The orders $\Ccal_{\omega_1}$ and $\Rcal_{\omega_1}$ force
$\diamond(\non(\Mcal))$ and
$\diamond(\non(\Ncal))$ respectively.
\end{thm}

\begin{proof}
The arguments for each are almost identical so we will only present the
case of $\Rcal_{\omega_1}$.
Let $\dot G$ be an $\Rcal_{\omega_1}$-name for the element of
$2^{\omega_1}$ corresponding to the generic filter.
Fix an $\Rcal_{\omega_1}$-name $\dot F$ for a Borel map
from $2^{<\omega_1}$ to $\Ncal$ and let
$\dot r_\delta$ be an $\Rcal_{\omega_1}$-name
for a real such that $\dot F \restriction 2^{\delta}$ is definable
from $\dot r_\delta$.
Pick a strictly increasing function $f:\omega_1 \to \omega_1$ such that
$\dot r_\delta$ is forced to be in $\Vbf[\dot G\restriction f(\delta)]$.
Let $\dot g(\delta)$ be defined to be
$\dot G \restriction [f(\delta),f(\delta) + \omega)$ (interpreted
canonically
as a real).

To see that $\dot g$ works, let $\dot f:\omega_1 \to 2$ be an
$\Rcal_{\omega_1}$-name.
Let $C$ be the collection of all $\delta$ for which it is forced
that $\dot f \restriction \delta \in \Vbf[\dot G \restriction \delta]$.
Because $\Rcal_{\omega_1}$ is c.c.c., $C$ is closed and unbounded.
Since $\dot G$ is generic, $\dot g(\delta)$ avoids every null set
coded in $\Vbf[\dot G \restriction f(\delta)]$, including
$\dot F(\dot f \restriction \delta)$.
\end{proof}

The above proof actually shows that $\diamond^*(\non(\Mcal))$ and
$\diamond^*(\non(\Ncal))$ hold in the corresponding models where
$\diamond^*(A,B,E)$ is obtained from $\diamond(A,B,E)$ by
replacing ``stationary'' by ``club.''
One could, of course, produce a myriad of results of a similar flavor:
e.g. $\diamond^*(\cof(\Mcal))$
holds after adding $\omega_1$ Hechler reals or $\diamond^*(\sfrak)$
holds
after generically adding a sequence of $\omega_1$ independent
reals.\par

It should be noted that the results of
\cite{thesis:Hirschorn}, \cite{random_top}, and \cite{random:Todorcevic}
place considerable limitations on the strength of
$\diamond^*(\non(\Ncal))$ --- and hence $\diamond(\non(\Ncal))$ ---
as they show that there are a number of consequences of $\MA_{\aleph_1}$
which are consistent with $\diamond^*(\non(\Ncal))$.
For instance Theorem \ref{Cohen_random_diamond} gives the
following corollary which contrasts
the remarks preceding Definition \ref{Borel_Tukey}.

\begin{thm}
\label{different_guessing}
It is relatively consistent with $\CH$ that $\diamond(\non(\Ncal))$
holds but $\diamond(\non(\Mcal))$ fails.
\end{thm}

\begin{proof}
By a result of Hirschorn \cite{thesis:Hirschorn}, it is consistent
with $\CH$ that after forcing with any measure algebra there are
no Suslin trees.
After forcing with $\Rcal_{\omega_1}$ over this model we obtain
a model in which $\diamondsuit(\non(\Mcal))$ fails but
$\diamondsuit(\non(\Ncal))$ holds.
\end{proof}

So in particular, $\diamond(\non(\Ncal))$ is not sufficient to
carry out the construction of a Suslin tree.

\begin{question}
Does $\diamond(\bfrak)$ imply the existence of a Suslin tree?
\end{question}

This also suggests the following meta-question:

\begin{question}
If $(A_1,B_1,E_1)$ and $(A_2,B_2,E_2)$ are two Borel invariants such
that the inequality $\<A_1,B_1,E_1\> < \<A_2,B_2,E_2\>$ is consistent,
is it consistent that
$\diamond(A_1,B_1,E_1)$ holds and $\diamond(A_2,B_2,E_2)$ fails in the
presence of $\CH$?
\end{question}

We will now move on to study countable support iterations.

\begin{defn}
A \emph{Borel} forcing notion is a partial order $(X,\leq)$
with a maximal element\footnote{In this paper we adopt the convention that
if $p$ is stronger than $q$ then we write $p \leq q$.}
such that $X$ and $\leq$ are Borel sets.
\end{defn}

Given a Borel forcing notion, we will always interpret it in forcing
extensions
using its code rather than taking the ground model forcing notion.
Observe that many Borel forcing notions $\Qcal$
designed for adding a single real
(e.g. those for adding Laver, Miller, Sacks, etc. reals) are
equivalent to the forcing $\Pcal(2)^{+} \times \Qcal$ where
$\Pcal(2)$ is considered as the Boolean algebra with two atoms.
The atom of $\Pcal(2)$ which the generic selects can be thought of
as the first coordinate of the generic real which is added.

The following theorem will now become our focus:

\begin{thm} \label{P_forces_diamonds}
Suppose that $\<\Qcal_\alpha:\alpha < \omega_2\>$
is a sequence of Borel partial orders such that for each $\alpha < \omega_2$
$\Qcal_\alpha$ is equivalent to $\Pcal(2)^+ \times \Qcal_\alpha$
as a forcing notion and let $\Pcal_{\omega_2}$ be the countable
support iteration of this sequence.
If $\Pcal_{\omega_2}$ is proper and $(A,B,E)$ is a Borel invariant then
$\Pcal_{\omega_2}$ forces $\<A,B,E\> \leq \omega_1$ iff
$\Pcal_{\omega_2}$ forces $\diamond(A,B,E)$.
\end{thm}

\begin{remark}
This is actually a rather weak formulation what can be proved.
All of ``Borel'' that is used is that the forcing notions
remain forcing notions in generic extensions and they can be computed
from a real.
Also, it is not entirely necessary that the forcing notions be in
$V$; we will need only that the choice of the sequence of
forcing notions does not depend on the ``first coordinates'' of
the first $\omega_1$ generic reals added by the iteration.
We chose the phrasing that we did both because of its simplicity and
the fact that it covers most countable support iterations of definable
forcings.
\end{remark}

We will prove this theorem as a series of lemmas.

\begin{lem} \label{1stCoord:lem}
Suppose that
$\<\Pcal_\alpha,(\Pcal(2)^+ \times \dot \Qcal_\alpha):\alpha < \omega_1\>$
is a countable support iteration such that for all $\alpha < \omega_1$
$\Pcal_\alpha$ forces $\dot \Qcal_\alpha$ is a proper partial order.
For all $\alpha \leq \omega_1$,
the suborder $\Pcal^0_\alpha \subseteq \Pcal_\alpha$ of all conditions
whose first coordinate is trivial is completely embedded in $\Pcal_\alpha$.
\end{lem}

\begin{proof}
By induction on $\alpha$ we prove that the identity map
$\iota_\alpha:\Pcal_\alpha^0 \to \Pcal_\alpha$ is a complete embedding.
Note that for $\gamma < \alpha$,
$\Pcal_\gamma^0 = \{p \restriction \gamma:p \in \Pcal^0_\alpha\}$.
Also observe that $\Pcal_\alpha^0$ the direct limit of
$\Pcal_\gamma^0$ $(\gamma < \alpha)$ under the usual system of embeddings.

If $\alpha = 0$ this is trivial.
By the above observation, for a limit ordinal $\alpha > 0$ we have
(checking conditions (1) -- (3) of Definition 7.1 of Ch VII in \cite{set_theory}):
\begin{enumerate}

\item Since $\iota_\alpha$ is the inclusion map, it automatically preserves order.

\item If $p,p'$ are incompatible in $\Pcal^0_\alpha$, there must be a $\gamma < \alpha$ such that
$p \restriction \gamma$ is incompatible with $p \restriction \gamma'$ in $\Pcal^0_\gamma$ (this is a standard fact about
direct limits --- see 5.11 of Ch VIII in \cite{set_theory}).
By the induction hypothesis, $p \restriction \gamma$ and $p' \restriction \gamma$ are incompatible in
$\Pcal_\gamma$ and hence $p$ and $p'$ are incompatible in $\Pcal_\alpha$.

\item Given $q$ in $\Pcal_\alpha$ let $(\dot \epsilon_\gamma,\dot q_\gamma)$ denote
$q(\gamma)$.
It is easy to check by induction on $\gamma < \alpha$ that there is a unique condition
$\bar q$ in $\Pcal^0_\alpha$ such that $\bar q(\gamma) = (\one,\dot q_\gamma)$
for all $\gamma < \alpha$.
We now claim that if $r \in \Pcal^0_\alpha$ extends $\bar q$ then $r$ is compatible with $q$.
Indeed, the $\gamma\Th$ coordinate of the common extension is
$(\epsilon_\gamma,\dot r_\gamma)$ where $r(\gamma) = (\one,\dot r_\gamma)$. 

\end{enumerate}
This finishes the limit case of the inductive proof; the successor case is similar.
\end{proof}

\begin{lem} \label{iterationToProd:lem}
Suppose $\Tcal$ is a forcing notion which does not
add any countable sequences of ordinals and that
$\Pcal = \<\Pcal_\alpha,\dot \Qcal_\alpha:\alpha < \delta\>$
is a proper countable support iteration of Borel forcing notions.
The forcings $\Tcal * \dot \Pcal$ and $\Pcal \times \Tcal$
are equivalent provided that they do not collapse $\omega_1$.
\end{lem}

\begin{proof}
The $\Tcal$-name $\dot \Pcal$ will be used to refer to the iteration computed
after forcing with $\Tcal$. 
It now suffices to show that $\dot \Pcal$ is equal to $\check \Pcal$.
This will be proved by induction on $\delta$.

If $\delta$ is a limit of uncountable cofinality
then, by the inductive hypothesis, $\check \Pcal$ and $\dot \Pcal$ are both
inverse limits of equal orders (computed before and after forcing with $\Tcal$
respectively).
Since $\Tcal$ adds no new countable sequences of ordinals,
it forces that $\cf(\delta) > \omega$.
Therefore the inverse limit construction is absolute and
we have that $\dot \Pcal$ equals $\check \Pcal$.

If $\delta$ is a limit of countable cofinality, the same argument applies with
the observation that direct limits of systems
with countable cofinality are absolute between models with the same countable
sequences of ordinals.

Finally, if $\delta = \alpha + 1$ then
$\Pcal = \Pcal_\alpha * \dot \Qcal_\alpha$.
First we will prove that $\Tcal$ does not
add any new $\Pcal_\alpha$-names for reals.
To this end, suppose that $t$ is in $\Tcal$, $\dot r$ is forced by $t$ to
be a $\Tcal$-name for a $\Pcal_\alpha$-name for an element of $2^\omega$.
Let $\dot A_n$ $(n < \omega)$ be a countable sequence
such that $t$ forces that $\dot A_n$ is a
maximal antichain in $\Pcal_\alpha$ whose elements
decide the value of $\dot r(n)$.
Since $\Tcal * \Pcal_\alpha$ does not collapse $\omega_1$,
there is a $\bar t$ extending
$t$, a $p$ in $\Pcal_\alpha$ and $\Tcal$-names $\dot C_n$
such that $\bar t$ forces
that $\dot C_n$ is a countable subset of $\dot A_n$ and is maximal below $p$.
Since $\Tcal$ does not add countable sequences,
there is an extension of $\bar t$
which decides $\dot C_n$ for all $n$.
Hence $\Tcal$ does not add any $\Pcal_\alpha$-names for reals and therefore
$\Qcal_\alpha$ is the same computed after forcing with
$\Tcal * \Pcal_\alpha$ as it is after forcing with $\Pcal_\alpha$.
Combining this with the inductive hypothesis we have
$$\check \Pcal = \check \Pcal_\alpha * \check \Qcal_\alpha =
\dot \Pcal_\alpha * \dot \Qcal_\alpha
= \dot \Pcal,$$
thus finishing the proof.
\end{proof}

\begin{defn}
A forcing notion $(\Pcal,\leq)$ is \emph{nowhere c.c.c.} if for every $p$ in $\Pcal$
there is an uncountable antichain of elements which extend $p$.
\end{defn}
 
\begin{lem} \label{T:lem}
If $\Pcal_{\omega_2}$ is as in the statement of Theorem \ref{P_forces_diamonds}
and $\Pcal_{\omega_1}$ is the c.s. iteration of $\<\Pcal_\alpha,\Qcal_\alpha:\alpha < \omega_1\>$ then
there is a $\Pcal_{\omega_1}$-name $\dot \Tcal$ for a tree of height $\omega_1$ which is
nowhere c.c.c. and does not add reals such that
$\Pcal_{\omega_2}$ is equivalent to $\Pcal_{\omega_2} * \dot \Tcal$.
\end{lem}

\begin{proof}
By passing to an equivalent iteration, we replace $\Pcal_{\omega_1}$
by the c.s. iteration of the orders $\Pcal(2)^+ \times \Qcal_\alpha$.
Let $\dot \Tcal$ be the $\Pcal^0_{\omega_1}$-name for the quotient of
$\Pcal_{\omega_1}$ by the $\Pcal^0_{\omega_1}$-generic filter.
Thus $\dot \Tcal$ is a tree of height $\omega_1$ in which the $\alpha\Th$ level
is the image of $\Pcal^d_\alpha$ in the quotient.
Observe that if $G \subseteq \Pcal^0_{\omega_1}$ is generic,
$\Tcal$ is the collection of all $t:\alpha \to 2$ in $V[G]$ such that
for every $\gamma \leq \alpha$, $t \restriction \gamma$ is in $V[G \cap \Pcal_\gamma^0]$.
From this it is clear that $\Tcal$ is nowhere c.c.c.. 
Since $\dot \Tcal$ is $\Pcal_{\omega_1}^0$-name for a tree of size $\omega_1$
which is everywhere of uncountable height and which
embeds into a proper partial order,
$\Pcal^0_{\omega_1}$ forces that $\dot \Tcal$ does not add any
new countable sequences of ordinals.

Let $\dot \Pcal'$ be the $\Pcal_{\omega_1}$-name for the
remaining part of the iteration $\Pcal_{\omega_2}$.
Now $\Pcal_{\omega_2}$ is equivalent to
$(\Pcal^0_{\omega_1} * \dot \Tcal) * \dot \Pcal'$ which is in turn
equivalent to $\Pcal_{\omega_1} * (\dot \Tcal * \Pcal')$ which is
equivalent to $(\Pcal_{\omega_1} * \dot \Pcal') * \dot \Tcal = 
\Pcal_{\omega_2} * \dot \Tcal$.
\end{proof}

The following lemma now completes the proof of Theorem \ref{P_forces_diamonds}.

\begin{lem}\label{organize}
Let $(A,B,E)$ be a Borel invariant such that $\<A,B,E\> \leq \omega_1$.
If $\Tcal$ is a tree of height $\omega_1$
which is nowhere c.c.c. and does not add reals then
$\Tcal$ forces $\diamond(A,B,E)$.
\end{lem}

\begin{proof}
Let $b_\xi$ $(\xi < \omega_1)$ be a sequence of elements of $B$ which
witnesses $\<A,B,E\> \leq \omega_1$ and let
$\dot F:2^{< \omega_1} \to A$ be a $\Tcal$-name for a Borel function.
For each $\delta < \omega_1$ pick a $\Tcal$-name $\dot r_\delta$ for
a real which codes $\dot F \restriction \delta$.
For each $t$ in $\Tcal$ of height $\delta$
pick a real $s_t$ and a map $h_t:\omega_1 \to \Tcal$ such that
\begin{enumerate}

\item the collection $\{h_t(\xi):\xi < \omega_1\}$ is an antichain and

\item $h_t(\xi)$ extends $t$ and forces $\dot r_\delta$ to be $s_t$.

\end{enumerate}
Define a $\Tcal$-name $\dot g$ for a function from $\omega_1$ into
$B$ by making $h_t(\xi)$ force that $\dot g(\delta) = b_\xi$
where $\delta$ is the height of $t$
(if $\dot g$ is undefined somewhere define it arbitrarily).

Now let $\dot f$ be a $\Tcal$-name for a function from $\omega_1$ to $2$
and $\dot C$ be a $\Tcal$-name for a closed unbounded subset of $\omega_1$.
Let $A_n$ be a sequence of maximal antichains in $\Tcal$ such that
if $u$ is in $A_n$ and has height $\delta$ and $\bar u$ is in $A_{n+1}$ and
extends $u$ then $\bar u$ decides $\dot f \restriction \delta$ and forces that
there is an element of $\dot C$ between $\delta$ and the height of $\bar u$.
Since $\Tcal$ does not add reals, there is a minimal $t$
such that for every $n$ there is a $u_n$ in $A_n$ which is below $t$.
Hence if $\delta$ is the height of $t$, $t$ decides $\dot f \restriction \delta$
and forces $\delta$ to be in $\dot C$.
Now there is an $a$ in $A$ such that $t$ forces that if
$\dot F(\dot f \restriction \delta)$ is computed using the code $s_t$ then
its value is $a$.
Find a $\xi$ such that $(a,b_\xi)$ is in $E$.
The condition $h_t(\xi)$ forces that $\delta$ is in $\dot C$ and that
$(\dot F(\dot f \restriction \delta), \dot g(\delta))$ is in $E$,
finishing the proof.
\end{proof}

The following is a typical corollary of the previous two theorems.
We will see in Section \ref{structured_invariants:sec} that this in turn
implies that $\afrak = \ufrak = \omega_1$ in the iterated Sacks model.

\begin{cor} \label{diamonds_in_Sacks}
Both $\diamond(\rfrak)$ and $\diamond(\dfrak)$ hold in the
iterated Sacks model.
\end{cor}
 
The following result gives another way of seeing the relative consistency
of $\clubsuit + \neg \CH$.\footnote{
Baumgartner has demonstrated in an unpublished note
that $\club$ holds in the Sacks model.
This result was obtained shortly after Shelah's proof
of the consistency of $\club + \neg \CH$ \cite{Baumgartner:email}.}
Unlike the standard proofs (see Chapter I Section 7 of
\cite{proper_forcing})
where one deliberately arranges that $\clubsuit$ holds in the forcing
extension, the Sacks model was considered for entirely different
reasons.

\begin{cor}
$\clubsuit$ holds in the iterated Sacks model.
\end{cor}

\begin{proof}
Without loss of generality we may assume that our ground model
satisfies $\CH$.
It suffices to show that $\Sbb_{\omega_2} * \dot \Tcal$ forces
$\clubsuit$ where $\Tcal$ is the forcing notion from Lemma \ref{T:lem}.
In \cite{life_in_Sacks}
it has been (essentially) shown that for every $\Sbb_{\omega_2}$-name
$\dot X$ for an uncountable subset of $\omega_1$ there is a
$\Sbb_{\omega_2}$-name $\dot C$ for a closed and unbounded subset of
$\omega_1$
such that if $p$ forces that $\delta$ is in $\dot C$ then there is a
$q$ extending $p$ and a ground model $A \subseteq \delta$ which is
cofinal such that $q$ forces that $A$ is contained in $\dot X$.

We will now work in the forcing extension given by $\Sbb_{\omega_2}$.
For each $t$ in $\Tcal$, let $h_t:\omega_1 \to \dot \Tcal$ be a 1-1
function such that the range of $h_t$ is an antichain above $t$.
For limit $\delta$ define a $\Tcal$-name $\dot C_\delta$ by letting
$h_t(\xi)$ force $\dot C_\delta = A_\xi$ where $\{A_\xi:\xi < \omega_1\}$
enumerates the cofinal subsets of $\delta$ before forcing with
$\Sbb_{\omega_2}$. 
The method of proof of Lemma \ref{organize}
now shows that $\dot C_\delta$ $(\delta \in \lim(\omega_1))$ is forced
to be a $\clubsuit$-sequence.
\end{proof}

One ``rule of thumb'' which one learns when working with the classical
invariants of the form $(A,B,E)$ is that,
if $\<A,B,E\> < \<C,D,F\>$ is consistent then this can typically
accomplished by a countable support
iteration of length $\omega_2$ of proper Borel
forcing notions in $\Vbf$ (typically the sequence $\Qcal_\alpha$
$(\alpha < \omega_2)$ is a constant sequence).\footnote{In general
this is a phenomenon which is not well understood and is currently being
analyzed by a number of people.
There are Borel invariants such as
$\cov(\Ncal)$ and $(\Rbb^\omega,\Ncal,\sqsubseteq)$ which can only be
separated if the continuum is larger than $\aleph_\omega$
($f \sqsubseteq E$ if the range of $f$ is contained in $E$).
This is because $\cov(\Ncal)$ can have countable cofinality \cite{cf(cov(N))} while
$(\Rbb^\omega,\Ncal,\sqsubseteq)$ cannot and yet
$\cov(\Ncal) \leq \<\Rbb^\omega,\Ncal,\sqsubseteq\>
\leq \cf([\cov(\Ncal)]^{\omega},\subseteq)$.
}
In such a case, Theorem \ref{P_forces_diamonds} tells
us that $\diamond(A,B,E)$ does not imply $\<C,D,F\>$ is $\omega_1$.
The reader is referred to \cite{set_theory_reals}
for an introduction to some of the common Borel forcing notions and
\cite{norms_possibilities} for some
of the more advanced techniques for building Borel forcing notions.

The above results imply that
$\diamond(\Rbb,\ne)$ holds
in many of the models obtained by adding a specific type of real.
The following theorem, however, gives a much more natural setting
for studying $\diamond(\Rbb,\ne)$ and its consequences.

\begin{thm}
After forcing with a Suslin tree $\diamond(\Rbb,\ne)$ holds.
\end{thm}

\begin{proof}
Similar to the proof of Theorem \ref{P_forces_diamonds}
(in fact it is most natural to show that $\diamond(\omega,=)$ holds
after forcing with a Suslin tree).
\end{proof}

Many of the combinatorial consequences of $2^{\omega} < 2^{\omega_1}$
are in fact consequences of $\diamond(\Rbb,\ne)$.
It should be noted that Farah, Larson, Todor\v{c}evi\'c
and others have noticed that these consequences
hold after forcing with a Suslin tree.

\begin{thm}
$\diamond(\Rbb,\ne)$ implies:
\begin{enumerate}

\item \label{t} $\tfrak = \omega_1$.

\item \label{Q} There are no $Q$-sets.

\item \label{ladder} Every ladder system has a non-uniformizable
coloring.

\item \label{K3} There is an uncountable subset of a c.c.c. partial
order
with no uncountable $3$-linked subcollection.

\end{enumerate}
\end{thm}

\begin{proof}
Item \ref{t} is deferred to Theorem \ref{t:thm} of the next section.
The proof that $\diamond(\Rbb,\ne)$ implies
items \ref{Q} and \ref{ladder} is the same
as the proof that $\Phi(2,\ne)$ implies these statements
(see \cite{weak_diamond}).
Item \ref{K3} can be extracted from the proof of Theorem 7.7 of
\cite{partition_problems} and Theorem \ref{t:thm} below.
\end{proof}

On the other hand, Larson and Todor\v{c}evi\'{c} have had a great deal
of success in proving that certain consequences of $\MA_{\aleph_1}$
and other forcing axioms can hold after forcing with a Suslin tree
(see \cite{katprob}, \cite{cc_max}).
A major open question in this line of research is:

\begin{question}
Is $\diamond(\Rbb,\ne)$ consistent with the assertion that
every c.c.c. forcing notion has Property K?
\end{question}

\section{$\diamond$-principles and cardinal invariants}

\label{structured_invariants:sec}

There are a number of well studied cardinal invariants of the continuum
which do not satisfy our definition of ``invariant.''
Generally this is because the invariants in question make reference to
some additional structure.
For instance, $\ufrak$ can be considered to be the smallest size
of a reaping family which is also a filter base.
A natural question to ask is how these cardinals are influenced by
the $\diamond$-principles we have considered thus far.
It turns out that these $\diamond$-principles do have a strong impact
on cardinals such as $\tfrak$, $\afrak$, and $\ufrak$ and moreover
provide a uniform approach for computing the values of these invariants
in many standard models.

The first instance of this influence was Hru\v{s}\'{a}k's proof that
$\diamond_{\dfrak}$ implies $\afrak = \omega_1$.
In addition to allowing for easier computations, the results
below explain why the proofs of statements such as
$\CON(\bfrak < \afrak)$ and $\CON(\rfrak < \ufrak)$
require more sophisticated arguments than, e.g.,
$\CON(\bfrak < \dfrak)$.
It also suggests that there are no natural formulations of statements
such as $\diamond(\tfrak)$ and $\diamond(\afrak)$.

The first theorem is essentially a recasting of the well known fact that
$2^\omega < 2^{\omega_1}$ implies $\tfrak = \omega_1$.

\begin{thm}\label{t:thm}
$\diamond(\Rbb,\ne)$ implies $\tfrak = \omega_1$.
\end{thm}

\begin{proof}
$\diamond(\Rbb,\ne)$ is equivalent to $\diamond(2,\ne)$ so
we will use this assumption instead.
Let $X$ be the subset of $([\omega]^\omega)^\omega$ consisting
of all strictly $\subseteq^*$ decreasing sequences of sets.
Let $D:X \to [\omega]^\omega$ be defined by setting
the $n\Th$ element of $D(\vec{A})$ to be the least element of
$\bigcap_{i \leq n} A_i$
which is greater than $n$.
Notice that $D(\vec{A})$ is almost contained in $A_n$ for all $n <
\omega$ and $D$ is continuous.

Our map $F$ will be defined on pairs $\vec{A},C$
where $\vec{A} = \<A_\xi: \xi < \delta\>$
is a strictly $\subseteq^*$-decreasing sequence,
$\delta$ is a limit and $C$ is an infinite subset of $\omega$
which is almost contained in $A_\xi$ for all
$\xi < \delta$.
Let $B(\vec{A})$ be the collection  of all even indexed
elements of $D(\<A_{\delta_n}:n \in \omega\>)$ in its increasing
enumeration.
Set $F(\vec{A},C)$ to be $0$ if $C$ is almost contained in $B(\vec{A})$ and
$1$ otherwise.

Let $g:\omega_1 \to 2$ be a $\diamond(\Rbb,\ne)$-sequence for $F$.
Construct $\<A_\xi:\xi \in \omega_1\>$ by recursion.
Let $A_n$ $(n \in \omega)$ be any strictly decreasing $\omega$-sequence in
$[\omega]^\omega$.
Now suppose that $\<A_\xi:\xi < \delta\>$ is given.
Define $A_\delta$ to be $B(\vec{A})$ if $g(\delta) = 0$ and
$D(\<A_{\delta_n}:n \in \omega\>) \setminus B(\vec{A})$
otherwise.
It is easily checked that if
$F(\<A_\xi:\xi < \delta\>,C)$ is defined and not equal to $g(\delta)$
then
$A_\delta$ does not almost contain $C$.
\end{proof}

The next result can be considered as an optimization of Theorem
\ref{diamond_d/a}.
It is an old result of Solomon that $\bfrak \leq \afrak$ is provable
in $\ZFC$ \cite{integers_topology}.

\begin{thm}
$\diamond(\bfrak)$ implies $\afrak = \omega_1$.
\end{thm}

\begin{remark}
Shelah has shown that $\bfrak < \afrak$ is consistent \cite{cardinal_invariants}
(see also \cite{mob_mad}).
\end{remark}

\begin{proof}
We will first define a Borel function $F$ into the set $\NN$ as follows.
The domain of $F$ is the set of all pairs $(\<A_\xi:\xi < \delta\>,B)$
such that:
\begin{enumerate}

\item $\delta$ is an infinite countable ordinal.

\item $\{A_\xi:\xi<\delta\} \cup \{B\}$
is an almost disjoint family of infinite subsets of $\omega$.

\item For infinitely many $n$ the set
$B \cap A_{e_\delta(n)} \setminus \bigcup_{i < n} A_{e_\delta(i)}$
is non-empty. 

\end{enumerate}
We will denote the set of $n$ from condition 3 by $I(\vec{A},B)$.
Define $$F(\<A_\xi:\xi < \delta\>,B)(k) =
\min\left(B \cap A_{e_\delta(n)} \setminus \bigcup_{i < n} A_{e_\delta(i)}\right)$$
where $n$ is the $k\Th$ least element of $I(\vec{A},B)$.

Now suppose that $g:\omega_1 \to \NN$
is a $\diamond(\bfrak)$-sequence for $F$.
By making the entries in $g$ larger if necessary,
we may assume that they form a $<^*$-strictly increasing
sequence of increasing functions.

We will now construct a maximal almost disjoint family by
recursion.
Let $\<A_n:n< \omega\>$ be any almost disjoint family of
infinite subsets of $\omega$.
If $\<A_\xi:\xi < \delta\>$ has been defined, set
$$A_\delta =
\omega \setminus \bigcup_{n < \omega}
\left[
A_{e_{\delta}(n)} \setminus \left(g(\delta)(n) \cup \bigcup_{i<n} A_{e_\delta(i)}
\right) \right].$$
Since for each $n$ the set $g(\delta)(n) \cup \bigcup_{i<n} A_{e_\delta(i)}$
has finite intersection with $A_{e_\delta(n)}$, $A_\xi$
has finite intersection with $A_\delta$ for each
$\xi < \delta$.

To see that $\{A_\xi:\xi < \omega_1\}$ is maximal, suppose that
$B$ is an infinite subset of $\omega$.
First notice that if $\delta$ is at least $\omega$
and $(\<A_\xi:\xi < \delta\>,B)$ satisfies condition 2 but not condition 3 then
$B$ has infinite intersection with (in fact is almost contained in)
$A_\delta$.
Therefore we will be finished if we can show that if
$(\<A_\xi:\xi <\delta\>,B)$ satisfies conditions 1-3 and
$g$ guesses $(\<A_\xi:\xi < \omega_1\>,B)$ at $\delta$ then
$B \cap A_\delta$ is infinite.

To this end, suppose
$$F(\<A_\xi:\xi <\delta\>,B) \not >^* g(\delta)$$
and let $N$ be a given natural number.
For ease of reading we will let $\vec{A}$ abbreviate
$\<A_\xi:\xi < \delta\>$.
Find a number $k$ such that the $k\Th$ least element $n$ of $I(\vec{A},B)$
has the following properties:
\begin{enumerate}

\item $g(\delta)(k)$ is greater than $F(\vec{A},B)(k)$

\item the minimum $l$ of $B \cap A_{e_\delta(n)} \setminus \bigcup_{i < n}
A_{e_\delta(i)}$ is greater then $N$.

\end{enumerate}
The last choice is possible since
$$\{A_{e_\delta(j)} \setminus \bigcup_{i < j}
A_{e_\delta(i)}:j < \omega\}$$ forms
a disjoint family of sets.
It is now sufficient to show that $l$ is in $A_\delta$.
Observe that the only possibility for removing $l$ from $A_\delta$
is with the index $n$ since $l$ is in every set of the form
$\bigcup_{i<m} A_{e_\delta(i)}$ for $m > n$ and not in any $A_{e_\delta(i)}$
for $i < n$.
Since $k \leq n$ and $g(\delta)$ is monotonic,
$l = F(\vec{A},B)(k) < g(\delta)(k) \leq g(\delta)(n)$.
Thus $l$ is not in
$$A_{e_{\delta}(n)} \setminus
\left( g(\delta)(n) \cup \bigcup_{i<n} A_{e_\delta(i)} \right)$$
and therefore is in $A_\delta$ as desired.
\end{proof}

\begin{notn}
If two functions $f,g$ in $\omega^\omega$ are equal infinitely
often, then we will write $f =^\infty g$.
\end{notn}

It is known that the cardinal $\non(\Mcal)$ is equal to
$\<\omega^\omega,=^\infty\>$ \cite{set_theory_reals}.
 
\begin{defn} \cite{class_mad}
The cardinal $\afrak_e$ is the smallest size of a
maximal collection $\Acal \subseteq \NN$ of eventually
different functions.
\end{defn}

It follows from the above remark that
$\afrak_e \geq \non(\Mcal)$ and it has been shown by
Brendle that strict inequality is consistent \cite{mad_iteration}.

\begin{thm}
$\diamond(\NN,=^\infty)$ implies $\afrak_e = \omega_1$.
\end{thm}

\begin{proof}
Let $A_n$ $(n \in \omega)$ be a fixed partition of $\omega$ into
infinite pieces.
The domain of $F$ will be all countable sequences $\<f_\xi:\xi < \delta\>$
of eventually different functions and
an $h \in \NN$ which is eventually different from every
$f_\xi$ $(\xi < \delta)$.
For convenience our $F$ will take values in $(\omega^2)^\omega$.
Set $F(\<f_\xi:\xi < \delta\>,h)(n)$ to be
$(k,h(k))$ where $k$ is the least integer in $A_n$ such that
$g(l) \ne f_{e_\delta(i)}(l)$ for all $i \leq n$ and $l \geq k$.

Let $g:\omega_1 \to (\omega^2)^\omega$ be a $\diamond(\NN,=^\infty)$-sequence
for $F$.
Construct a sequence of eventually different functions
$f_\xi$ $(\xi \in \omega_1)$ by recursion.
Let $f_\xi$ for $\xi < \delta$ be a given sequence of eventually
different
functions.
Let $\Gamma$ be the collection of all $(k,v)$ such that
$k$ is in $A_n$, $g(\delta)(n) = (k,v)$, and
if $\xi < \delta$ with $e^{-1}_\delta(\xi) \leq n$ then
$f_\xi(k) \ne v$.
Notice that for a given $k$ there is at most one $v$ such that
$(k,v)$ is in $\Gamma$, and that $\Gamma$ is almost disjoint from
$f_\xi$ for all $\xi < \delta$.
Define $f_\delta(k)$ to be $v$ if $(k,v)$ is in $\Gamma$ for some $v$
and
$f_\delta(k)$ to be the least integer greater than
$f_\xi(k)$ for all $\xi$ with $e^{-1}_\delta(\xi) \leq k$.
Notice that $f_\delta$ is eventually different from $f_\xi$ for all
$\xi < \delta$.
To see that $\{f_\xi:\xi \in \omega_1\}$ is maximal, let $h \in
\omega^\omega$
and notice that
if $F(\<f_\xi:\xi < \delta\>,h)$ is defined and infinitely often equal to 
$g(\delta)$ then $f_\delta$ agrees with $h$ on an infinite set
--- namely those $k$'s for which $\Gamma$ was used in the definition of
$f_\delta(k)$.
\end{proof}

Recall that $\diamond_{\dfrak}$ is the following statement from
\cite{another_diamond}:
\begin{description}

\item[$\diamond_{\dfrak}$] There is a sequence $g_\delta:\delta \to
\omega$
indexed by $\omega_1$ such that if $f:\omega_1 \to \omega$ then there
is a $\delta \geq \omega$ such that
$f \restriction \delta <^* g_\delta$.

\end{description}
It is straightforward to check that $\diamond_{\dfrak}$ is a consequence of
$\diamond(\dfrak)$.
The following theorem answers a question asked in \cite{another_diamond}.

\begin{thm}
$\diamond_\dfrak$ implies that $\omega^\omega$ can be partitioned
into $\omega_1$ compact sets.
\end{thm}

\begin{remark}
Spinas has shown that it is consistent that $\dfrak = \omega_1$ and
yet $\omega^\omega$ cannot be partitioned into $\omega_1$
disjoint compact sets \cite{partition_numbers}.
\end{remark}

\begin{proof}
Notice first that any $\sigma$-compact subset of $\omega^\omega$ can
be partitioned into countably many compact sets.
This follows from the fact that $\omega^\omega$ is 0-dimensional.
If $f \in \NN$, let $K_f$ be the collection of all $g$ in $\NN$ such
that
$g \leq f$.
If $C \subseteq \NN$ is compact and $f \in \NN \setminus C$,
let $\Delta(f,C)$ be the maximum of $\Delta(f,y)$ where $y$ ranges
over $C$ (if $C$ is empty then let $\Delta(f,C) = 0$).
Since $C$ is compact and $f$ is not in $C$, this is always a finite
number.

Let $g_\delta$ $(\delta \in \omega_1)$ be a $\diamond_{\dfrak}$-sequence.
Given $C_\xi$ $(\xi < \delta)$, a disjoint sequence of compact sets for
limit
$\delta$, define
$$F_\delta = \bigcup_{g =^* g_\delta}
[K_{g \restriction \omega} \setminus \bigcup_{\xi < \delta}
\{x \in \NN:\Delta(x,C_\xi) > g(\xi)\}].$$
Notice that $\bigcup_{\xi < \delta} \{x \in \NN:\Delta(x,C_\xi) > g(\xi)\}$
is open and hence $F_\delta$ is $\sigma$-compact.
Let $\{C_{\delta+n}:n \in \omega\}$ be a partition of $F_\delta$ into
disjoint compact sets.
Clearly the sequence $C_\xi$ $(\xi \in \omega_1)$ is pairwise disjoint.
Let $x$ be in $\NN$ and suppose that $x$ is not contained in $C_\xi$ for
any
$\xi \in \omega_1$.
Define $f:\omega_1 \to \omega$ by setting $f \restriction \omega = x$
and
$f(\xi) = \Delta(x,C_\xi)$ if $\xi \geq \omega$.
Now pick an $\delta > \omega$ such that $f \restriction \delta <^* g_\delta$.
It follows that $x$ is in $F_\delta$ and therefore in $C_{\delta + n}$
for
some $n$, a contradiction.
\end{proof}

Recall that a free ultrafilter $\Ucal$ on $\omega$ is a \emph{P-point} if
whenever $F_n$ $(n \in \omega)$ is a sequence of elements of $\Ucal$, 
there is a $U$ in $\Ucal$ such that
$U\setminus U_n$ is finite for each $n\in\omega$.

\begin{thm}
$\diamond(\rfrak)$ implies that there is a P-point of character
$\omega_1$.
In particular $\diamond(\rfrak)$ implies $\ufrak = \omega_1$.
\end{thm}

\begin{remark}
Shelah and Goldstern have shown that $\omega_1 = \rfrak < \ufrak$ is
consistent \cite{r<u}.
\end{remark}

\begin{proof}
The domain of the function $F$ we will consider will consist of
pairs $(\vec{U},C)$ such that $\vec{U} = \<U_\xi:\xi < \delta\>$
is a countable $\subseteq^*$-decreasing
sequence of infinite subsets of $\omega$ and $C$ is a subset of $\omega$.

Given $\vec{U}$ as above, let
$B(\vec{U})$ be the set $\{k_i:i \in \omega\}$ where
$$k_i=\min(\bigcap_{j \leq i} U_{e^{-1}_\delta(j)}\setminus (k_{i-1}+1)).$$
Note that $B(\vec{U})$ is infinite and
almost contained in $U_\xi$ for every $\xi < \delta$.
Let $$F(\vec{U},C)=\{i: k_i\in C \cap B(\vec{U})\}$$
if  $\{i: k_i\in C \cap B(\vec{U})\}$
is infinite and let $$F(\vec{U},C)=\{i: k_i\not\in C \cap B(\vec{U})\}$$
otherwise.
Now suppose that $g:\omega_1 \to [\omega]^\omega$ is a
$\diamond(\rfrak)$-sequence for $F$.
Construct a $\subseteq^*$-decreasing sequence
$\<U_\xi: \xi \in \omega_1\>$ of
infinite sets by recursion. Let $U_n=\omega\setminus n$.
Having defined $\vec{U} = \<U_\xi:\xi < \delta\>$
let $U_\delta = \{ k_i: i\in g(\delta)\}$ where  $B(\vec{U})=\{k_i:i \in \omega\}$.

The family $\<U_\xi: \xi \in \omega_1\>$ obviously generates a P-filter. To see that
it is an ultrafilter, 
note that if a $C \subseteq \omega$ is  given and $g$ guesses $\vec{U},C$ at $\delta$
then $U_\delta$ is either almost contained in or almost disjoint from
$C$.
\end{proof}

By combining the  above proof with the argument that shows that 
$\dfrak=\omega_1$ implies the existence of a Q-point
one can without much difficulty prove the following.

\begin{cor}
$\diamond(\rfrak) + \dfrak=\omega_1$ implies that there is a selective
ultrafilter of character $\omega_1$.
\end{cor}

Recall that $\ifrak$ is the smallest cardinality of a maximal independent
family.
In \cite{balhr3} the rational reaping number
$$\rfrak_{\Qbb} = (\Pcal(\Qbb) \setminus \mathrm{NWD},\textrm{``does not
reap''})$$
is considered and it is proved that
$\rfrak,\dfrak \leq \rfrak_{\Qbb} \leq \ifrak$.
As with the earlier lemmas we show that the last inequality is,
in a sense, sharp.

\begin{thm}
$\diamond(\rfrak_{\Qbb})$ implies $\ifrak = \omega_1$.
\end{thm}

\begin{proof}
For this proof we will view $\Qbb \subseteq 2^\omega$ as the collection of all
binary sequences with finite support.
We will now define a Borel function $F$ on pairs
$(\<I_\xi:\xi < \delta\>,A)$ where
$\delta$ is an ordinal less than $\omega_1$ and $A$ and $I_\xi$ are subsets
of $\omega$ for all $\xi < \delta$.
The range of $F$ will be contained in $\Pcal(\Qbb)$.

If $\delta$ is finite or $\vec{I} = \<I_\xi:\xi < \delta\>$
is not independent then
return $\Qbb$ as the value of $F(\vec{I},A)$.
Otherwise, let $x_n(\vec{I})$ be the element of $2^\omega$ 
defined by $x_n(\vec{I})(k) = 1$ iff $n$ is in $I_{e_\delta(k)}$.
Observe that $X(\vec{I}) = \{x_n(\vec{I})\}_{n=0}^\infty$ is dense in $2^{\omega}$
since $\vec{I}$ is independent.
Fix a recursive homeomorphism $h$ from $X(\vec{I})$ to $\Qbb$.
Now put $F(\vec{I},A)$ to be the image of
$\{x_n(\vec{I}):n \in A\}$ under the map $h$.

Now suppose that $g$ is a $\diamond(\rfrak_{\Qbb})$-sequence for $F$.
We will now build an independent family $\{I_\xi :\xi < \omega_1\}$ by recursion.
Let $\{I_n :n < \omega\}$ be any countable independent family.
Now given $\vec{I} = \<I_\xi:\xi < \delta\>$, let
$t$ in $2^{< \omega}$ be such that $g(\delta)$ is dense in
$[t] = \{x \in 2^\omega:t \subseteq x\}$.
By altering $g(\delta)$ if necessary, we may assume that $h^{-1}(g(\delta))$ is contained
in $[t]$ and that $[t] \setminus h^{1}(g(\delta))$ is also dense in $[t]$.
Let $C = \{n \in \omega:h(x_n) \in g(\delta)\}$.
First we will see that $C$ has a nonempty intersection with
$\bigcap_{i < |u|} I^{u(i)}_{e_\delta(i)}$ iff $u$ extends $t$ where
$I^1 = I$ and $I^0 = \omega \setminus I$.
If $n$ is in such an intersection then $x_n(\vec{I})$ must be in $[u]$ by definition.
Now, if $u$ extends $t$, pick an $n$ such that $x_n$ is in $h^{-1}(g(\delta)) \cap [u]$.
Then $n$ is in $C$ and in $\bigcap_{i < |u|} I^{u(i)}_{e_\delta(i)}$.
Similarly one shows that $\omega \setminus C$ intersects every set of
the form $\bigcap_{i < |u|} I^{u(i)}_{e_\delta(i)}$ for every $u$ in $2^{< \omega}$.
Form $I_\delta$ so that $\<I_\xi:\xi \leq \delta\>$ is independent
and $I_\delta \cap \bigcap_{i < |t|} I^{t(i)}_{e_\delta(i)} = C$.

We are now finished once we show that $\{I_\xi:\xi < \omega_1\}$ is a maximal
independent family.
It is now sufficient to show that if $g$ guesses
$(\<I_\xi:\xi< \omega_1\>,A)$ at
$\delta$ then $\{I_\xi :\xi \leq \delta\} \cup \{A\}$ is not independent.
In fact, if $t$ is the element of $2^{<\omega}$
used in the definition of $C$ then
$$I_\delta \cap \bigcap_{i < |t|} I^{t(i)}_{e_\delta(i)}$$
is either contained in or disjoint from $A$.
\end{proof}

A natural question to ask is
whether $\diamond(\non(\Mcal))$ implies the existence of a Luzin set.
The answer is negative as the following theorem shows.

\begin{thm}
It is consistent that both $\diamond(\non(\Mcal))$
and $\diamond(\non(\Ncal))$ hold and there are no Luzin or
Sierpi\' nski sets.
\end{thm}

\begin{proof}
In the Miller model the cardinals
$\non(\Mcal)$ and $\non(\Ncal)$ are both $\omega_1$ and
hence by Theorem \ref{P_forces_diamonds} the corresponding
$\diamond$-principles hold.
On the other hand, Judah and Shelah \cite{killing_Luzin_Sierpinski}
have shown there are no Luzin or Sierpi\' nski sets in this model.
\end{proof}

One should note that results of this section combined with those of the
previous section
provide unified approach to determining values of cardinal invariants
with structure in 
many models in which this was traditionally done by arguments specific
to the forcing
construction at hand (see e.g. \cite{more_settheory_top}).

\section{$\diamond(A,B,E)$ and the Continuum Hypothesis}
\label{nonstandard_invariants:sec}

One of the most remarkable facts about the principle $\Phi$ of
Devlin and Shelah is that, while it resembles
a guessing principle in its statement, it is in fact equivalent to
the inequality $2^\omega < 2^{\omega_1}$ \cite{weak_diamond}.
The purpose of this section is to show that this phenomenon
is rather unique to the invariants
between $(\Rbb,\ne)$ and $(2,\ne)$ which characterize $\Phi$.
In particular we will show that $\diamond(\Rbb^\omega, \not \sqsupseteq)$
is not a consequence of $\CH$.
To emphasize the relevance of this to the invariants considered in
literature, we introduce the following definition.

\begin{defn} A Borel invariant $(A,B,E)$ is a
\emph{$\sigma$-invariant} if it satisfies the following strengthenings
\begin{description}

\item[$3^+.$] There is a Borel map $\Delta:A^\omega \to B$ such that
for all $\{a_n\}$ in $A^\omega$ the relation
$a_n E \Delta(\{a_n\})$ holds for all $n$.

\item[$4^+.$] There is a Borel map $\Delta^*:B^\omega \to A$ such that
for all $\{b_n\}$ in $B^\omega$ the relation
$\Delta^*(\{b_n\}) E b_n$ \emph{does not} hold for all $n$.

\end{description}
\end{defn}

Notice that if $(A,B,E)$ is a Borel $\sigma$-invariant then
$\<A,B,E\> \geq \omega_1$.
Moreover the cardinal invariants $(A,B,E)$ which appear in the literature
typically satisfy these conditions.
The connection to $(\Rbb^\omega,\not \sqsupseteq)$ is the following.

\begin{prop}
If $(A,B,E)$ is a Borel $\sigma$-invariant then $(A,B,E)$
is above $(\Rbb^\omega,\not \sqsupseteq)$ and below
$(\Rbb^\omega,\sqsubseteq)$ in the Borel Tukey order.
\end{prop}

\begin{proof}
We will only present the proof for $(\Rbb^\omega,\not \sqsupseteq)$
as the proof dualizes to the other case.
Since $B$ must be an uncountable Borel set for $(A,B,E)$ to
be a $\sigma$-invariant, we can find a Borel isomorphism between
$B$ and $\Rbb$.
Thus it suffices to show that
$(B^\omega,\not \sqsupseteq)$ is below $(A,B,E)$ in the Borel
Tukey order.
The map $f:B^\omega \to A$ is the map $\Delta^*$ and the map
$g:B \to B^\omega$ sends $b$ to the constant sequence $\bar b$.
If $\vec{b}$ is in $B^\omega$ and $b$ is in $B$ with
$\Delta^*(\vec{b}) E b$ then $b$ could not be in the range of $\vec{b}$.
Hence $\vec{b} \not \sqsupseteq \bar b$. 
\end{proof}

\begin{thm}
\label{guessing}
$\CH$ does not imply
$\diamond(\Rbb^\omega,\not \sqsupseteq)$.
\end{thm}

In order to prove this theorem, we will prove a more
technical result which may be of independent interest.

\begin{thm} \label{decomp_Pcr}
It is relatively consistent with $\CH$ that whenever
$\vec{C}$ is a ladder system and
$\vec{r}$ is a sequence
of distinct elements of $2^\omega$ indexed by $\lim(\omega_1)$, there
is a countable decomposition $\lim(\omega_1) = \bigcup_{n=0}^\infty X_n$
such that if $\gamma < \delta$ are in $X_n$ then
$$|C_\gamma\cap C_\delta| \leq \Delta(r_\gamma,r_\delta)$$
where $\Delta(r_\gamma,r_\delta)$ is the size of the largest common
initial segment of $r_\gamma$ and $r_\delta$.
\end{thm}

We shall now prove several lemmas that will at the end allow us to prove
Theorem \ref{decomp_Pcr}.
For the purposes of this proof let $\theta$ be a large enough regular
cardinal.

Given a pair $\vec{C}, \vec{r}$ 
as in the theorem we will denote by $\Qcal_{\vec{C},\vec{r}}$ the partial
order consisting of functions
$q:\lim(\delta +1)\rightarrow 2^{<\omega}$ for $\delta<\omega_1$ and
such that

\begin{enumerate}

\item $q(\beta)$ is an initial segment of $r_\beta$ for every
$\beta\in\dom (q)$,
 
\item if $\beta$ and $\gamma$ are distinct and
$q(\beta)=q(\gamma) $ then $| C_\gamma\cap C_\beta| \leq
\Delta(r_\gamma,r_\beta)$.

\end{enumerate}
ordered by extension (reverse inclusion). 

If $q$ is in $\Qcal_{\vec{C},\vec{r}}$,
$C\subseteq\omega_1$ has order type $\omega$,
$r$ is in $2^\omega$, and $\sigma$ is in $2^{<\omega}$,
then we will say that \emph{$q$ is consistent with
$(C,r)\mapsto \sigma$} if $\sigma$ is an initial segment of $r$ and 
$|C\cap C_\delta | \leq \Delta(r,r_\delta)$ 
for every $\delta\in \dom(q)$ with $q(\delta)=\sigma$.

\begin{lem} \label{basic_Pcr} Let $q \in \Qcal_{\vec{C},\vec{r}}$.
\begin{enumerate} 

\item If $C\subseteq\omega_1$ has order type $\omega$, $|C\cap
\dom(q)|$ is finite, and $r$ is in $2^\omega$ then
there is an $n \in \omega$ such that for all
$m\geq n$ $q$ is consistent with $(C,r) \mapsto r \restriction m$.

\item If $q$ is consistent with $(C_\alpha,r_\alpha)\mapsto \sigma$
then there is a $\bar q \leq q$ such that $\alpha$ is in the domain of
$\bar q$ and $\bar q(\alpha) = \sigma$.

\item Let $M$ be a countable elementary submodel of $H(\theta)$ such that
$\vec{r}\in M$.
If $\delta = M \cap \omega_1$ then for every 
$\beta \in M \cap \omega_1$ and $n \in \omega$
there is a $\gamma > \beta$ in $M$ such 
that $r_\gamma \restriction n =r_\delta \restriction n$.

\end{enumerate}
\end{lem}

\begin{proof} For (1) $n=|C\cap \dom(q)|$ obviously works. For (2)
enumerate $\lim(\alpha+1)\setminus\dom(q)$ as $\{\alpha_i: i\in I\}$,
where $I$ is either an integer or $\omega$, so that $\alpha=\alpha_0$.
Recursively pick a sequence
$\sigma_n$ ($n\in I$) so that 
$\sigma_0=\sigma$, $q$ is consistent with
$(C_{\alpha_n},r_{\alpha_n})\mapsto \sigma_n$, and 
$|\sigma_n|<|\sigma_{n+1}|$.
Then set
$$\bar q(\beta) =\begin{cases}
q(\beta) & \text{ if } \beta\in\dom(q)\\
                                  \sigma_n & \text{ if } \beta=\alpha_n
\text{ for some } n\in I.\end{cases}$$
To prove (3), let $\sigma=r_\delta\restriction n$. As it is finite,
$\sigma\in M$, and 
$\delta$ witnesses that
\[
H(\theta)\models(\exists\gamma>\beta)r_\gamma\restriction n = \sigma.
\]
Hence $M$ satisfies the same by elementarity.
\end{proof}

Notice first that if  $G$ is a $\Qcal_{\vec{C},\vec{r}}$-generic
filter then
$\{X_\sigma:\sigma\in 2^{<\omega}\}$ is the required decomposition,
where 
$X_\sigma= \{\alpha: \exists q\in G (q(\alpha)=\sigma)\}$.

Next we will show that the forcing $\Qcal_{\vec{C},\vec{r}}$ is
proper and does not add new reals,
and that moreover these forcings can be iterated
with countable support without adding reals.
Recall that a forcing notion $\Qcal$ is \emph{totally proper}
if for every countable elementary
submodel $M$ of $H(\theta)$ such that
$\Qcal \in  M$ and for every $q \in M \cap \Qcal$ there is a
$\bar q\leq q$ which is a lower bound for a filter containing $q$ which is
$\Qcal$-generic over $M$. 
Every such  $\bar q$ is called \emph{totally
$(M,\Pcal)$-generic}.
$\Qcal$ is  \emph{$\alpha$-proper} ($\alpha<\omega_1$) if for every
$q$ in $\Qcal$ and every increasing $\in$-chain
$\{M_{\beta}:\beta\leq\alpha\}$ of elementary
submodels of $H(\theta)$ such that $q,\Qcal\in M_0$,
there is a $\bar q\leq q$ which is
$(M_{\gamma},\Qcal)$-generic for every $\gamma\leq\alpha$. 
If $\Qcal$ is $\alpha$-proper for every $\alpha<\omega_1$ we will say
that
$\Qcal $ is\emph{ $<\omega_1$-proper}.\par

It is not difficult to see that a forcing notion $\Qcal$
is totally proper if and only if it is proper and does not add reals
(see \cite{no_ostaszewski}).

\begin{lem}\label{totally_proper_Pcr} The forcing notion
$\Qcal_{\vec{C},\vec{r}}$ is totally proper.
\end{lem}

\begin{proof} Let $M$ be an elementary submodel of $H(\theta)$ such that
$\vec{C},\vec{r}\in M$.
Fix a $q\in\Qcal_{\vec{C},\vec{r}}\cap M$ and 
select an enumeration $\{D_n:n\in\omega\}$ of all dense open subsets of
$\Qcal_{\vec{C},\vec{r}}$ which are elements of $M$.
Without loss of generality $M$ is an increasing
union of an $\in$-chain of elementary submodels $M_n$ $(n \in \omega)$
such that
$q$ and $\Qcal_{\vec{C},\vec{r}}$ are in $M_0$ and $D_n$ is in $M_n$.
Set $\delta = M \cap \omega_1$ and 
$\sigma = r_\delta \restriction |C_\delta\cap \dom(q)|$.
Construct a sequence $q_n$ $(n \in \omega)$ of conditions
together with a sequence $\beta_n$ $(n \in \omega)$ of 
ordinals so that
\begin{enumerate}

\item $q \geq q_0 \geq q_1 \geq \dots \geq q_n \geq q_{n+1} \geq \dots$,

\item $C_\delta \cap M_n \subseteq \beta_n \in M_n$,

\item $\Delta(r_{\beta_n},r_\delta) \geq |C_\delta \cap M_n|$,

\item $q_n \in M_n \cap D_n$, $\dom(q_n) = \beta_n+1$, and
$q_n(\beta_n)=\sigma$.

\end{enumerate}

Constructing these sequences is straightforward using Lemma \ref{basic_Pcr}.
Having done this let
$$\bar q=\bigcup_{n\in\omega} q_n\cup\{ (\delta,\sigma)\}.$$
Notice that as $q_n$ is consistent with $(C_\delta, r_\delta)\mapsto
\sigma$ for every $n\in \omega$, 
$\bar q$ is a condition in $\Qcal_{\vec{C},\vec{r}}$.
$\bar q$ is obviously totally
$(M,\Qcal_{\vec{C},\vec{r}})$-generic 
since for every $n$, $\bar q$ is below $q_n$ which is in $D_n$.

\end{proof}

\begin{lem}\label{<w_1-proper_Pcr} $\Qcal_{\vec{C},\vec{r}}$ is 
$\alpha$-proper for every $\alpha<\omega_1$.
\end{lem}

\begin{proof} We will prove the lemma by induction on $\alpha$.
Assume that the lemma holds for every $\beta < \alpha$.
Let $\{M_{\beta} : \beta \leq \alpha\}$ and
$q \in \Qcal_{\vec{C},\vec{r}} \cap M_0$ be given.
Set $\delta_{\beta} = M_{\beta} \cap \omega_1$
for each $\beta\ge\alpha$ and let $\sigma = r_{\delta_\alpha} \restriction
|C_{\delta_\alpha}\cap \dom(q)|$.

If $\alpha = \beta + 1$ for some
$\beta$ let $q' \in M_{\alpha}$, $q' \leq q$, be generic over all
$M_{\gamma}$ with $\gamma \leq \beta$.
As in Lemma \ref{totally_proper_Pcr} extend $q'$ to $\bar q$
which is (totally) generic over $M_{\alpha}$.

If $\alpha$ is a limit ordinal,
we will mimic the proof of Lemma \ref{totally_proper_Pcr}.
Fix a sequence of ordinals
$\alpha_n$ $(n \in \omega)$ increasing to $\alpha$.
Let $\{D_n : n \in \omega\}$ be an enumeration of all dense open
subsets of $\Qcal_{\vec{C},\vec{r}}$ in $M_\alpha$ such that $D_n\in
M_{\alpha_n}$.
Construct a sequence $q_n$ ($n\in\omega$) of conditions together with a
sequence $\beta_n$ $(n \in \omega)$ of 
ordinals so that
\begin{enumerate}

\item $q \geq q_0 \geq q_1 \geq \dots \geq q_n \geq q_{n+1} \geq \dots$,

\item $C_{\delta_\alpha}\cap M_{\alpha_n} \subseteq \beta_n \in
M_{\alpha_n}$,

\item $\Delta(r_{\alpha_n},r_{\delta_\alpha})\geq |C_{\delta_\alpha}\cap
M_{\alpha_n}|$,

\item $q_n \in M_n \cap D_n$, $\dom(q_n) = \beta_n + 1$, and
$q_n(\beta_n) = \sigma$, and

\item $q_{n+1}$ is $M_\gamma$-generic for every $\gamma \leq \alpha_n$.

\end{enumerate}
Let $\bar q=\bigcup_{n \in \omega} q_n \cup \{ (\delta_\alpha,\sigma)\}$. The
verification that
this works is as in Lemma \ref{totally_proper_Pcr}. 

\end{proof}

Recall the following definition and theorem from \cite{no_ostaszewski}.

\begin{defn} \label{2-complete} \cite{no_ostaszewski}
Let $\Pcal$ be totally proper and
$\dot\Qcal$ a $\Pcal$-name for a forcing notion and let $\theta$ be a
large enough regular cardinal. 
We shall say
that  $\dot\Qcal$ is \emph{$2$-complete for $\Pcal$} if {WHENEVER}
\begin{enumerate}

\item $N_0\in N_1 \in N_2$ are countable elementary submodels of
$H(\theta)$,

\item $\Pcal, \dot\Qcal\in N_0$,

\item $G\in N_1$ is $\Pcal$-generic over $N_0$ and has a lower
bound, and

\item $\dot q\in N_0$ is a $\Pcal$-name for a condition in
$\dot\Qcal$,

\end{enumerate}
{IT FOLLOWS THAT}
there is a $G'\in \Vbf$ which is $\dot\Qcal$-generic over
$N_0[G]$ such that $\dot q[G]\in G'$ and if $t\in \Pcal$ is a lower
bound for $G$ and $t$ is $\Pcal$-generic for $N_1$ and $N_2$ then
$t$ forces that $G'$ has a lower bound in $\dot \Qcal$.
\end{defn}

\begin{thm} \label{iter_tod} \cite{no_ostaszewski}
Let $\Pcal_\kappa = \<\Pcal_\alpha,\dot\Qcal_\alpha : \alpha < \kappa\>$ be a
countable support iteration such that
$\Vdash_\alpha$ ``$\dot\Qcal_\alpha$ is $<\omega_1$-proper and
$\dot\Qcal_\alpha$ is $2$-complete for 
$\Pcal_\alpha$". Then $\Pcal_\kappa$ is totally proper.
\end{thm}

\begin{lem}\label{2-complete_Pcr} Let $\Pcal$ be a totally proper
$<\omega_1$-proper poset and let $\dot\Qcal$ be a $\Pcal$-name for
$\Qcal_{\vec{C},\vec{r}}$ for some pair $\vec{C},\vec{r}$.
Then $\dot\Qcal$ is $2$-complete for $\Pcal$. 
\end{lem}

\begin{proof} Let $N_0 \in N_1 \in N_2$ be countable elementary submodels of
$H(\theta)$ and let  $\Pcal, \dot \Qcal \in N_0$.
Assume that $G \in N_1$ is an $(N_0, \Pcal)$-generic filter having a lower
bound and let $\dot q \in N_0$ be a $\Pcal$-name for a condition in
$\dot \Qcal$.
We have to find a $G'$ which is a $\dot\Qcal[G]$-generic
filter over $N_0[G]$ such that whenever $t \in \Pcal$ is a lower
bound for $G$ which is also $\Pcal$-generic over $N_1$ and $N_2$
then there is a $\Pcal$-name $\dot s$ such that $t \Vdash `` \dot
s$ is a lower bound for $G'$.''

Let $\delta = \omega_1 \cap N_0$ and set
$$\Dcal =
\{D\in N_0[G]: N_0[G]\models
\textrm{``} D \textrm{ is dense open in } \dot \Qcal [G] \textrm{''}\}.$$
Since $N_0,\dot Q,\Dcal$ and $G$ are all elements of $N_1$
and  $N_1\models\textrm{``} \Dcal
\textrm{ is countable''}$,
we can find an enumeration $\Dcal =\{ D_n: n\in\omega\}$ which is in $N_1$.

Let $\Ecal$ be the collection of all
$(C,r)$ such that $C$ is a cofinal subset of $\delta$ of order type
$\omega$ and $r$ is in $2^\omega$.
Clearly $\Ecal$ is in $N_1$.
Find an enumeration $\{(C^n,r^n):n \in \omega\}=$
$$\{(C,r)\in\Ecal\cap
N_1:(\forall \beta<\delta)(\forall n\in\omega)
(\exists \gamma\in [\beta,\delta))(r_\gamma\restriction n =r\restriction
n)\}$$
which is in $N_2$.
If we knew what $\dot C_\delta$ and
$\dot r_\delta$ evaluated to, we could proceed as in the proof of Lemma
\ref{totally_proper_Pcr} to produce $G'$. 
This is typically not the case.
What we do know, however, is that any $t$ which is a lower bound for $G$
and is $\Pcal$-generic over $N_1$ and $N_2$ forces that
$(\dot C_\delta, \dot r_\delta)$
appears in the enumeration $\{(C^n,r^n):n\in\omega\}$,
since $\Pcal$ does not add any new reals.
This allows us to simulate the proof of Lemma \ref{totally_proper_Pcr}
by diagonalizing over all possible
choices of $(\dot C_\delta, \dot r_\delta)$.

Again we may and will
assume that $N_0$ is the union of an $\in$-chain of elementary
submodels $M_n$ $(n \in \omega)$ such that
$\{M_n : n \in \omega\}$ is in $N_1$,
$\dot q[G]$ is in $M_0$ and $D_n$, $(C^n,r^n)$ are both in $M_n[G]$.
Construct a sequence $q_n$ $(n \in \omega)$ of conditions together with a
sequence $F_n$ $(n \in \omega)$ of 
finite sets of ordinals and $\sigma_n$ $(n \in \omega)$ of elements
of $2^{<\omega}$ by recursion on $n$ so that for every $i \leq n$
\begin{enumerate}

\item $\dot q[G] \geq q_i \geq q_n$,

\item  $C^i \cap M_n[G] \subseteq \min F_n \in M_n[G]$,
\item  $q_n$ is consistent with $(C^i,r^i) \mapsto \sigma_i$
\item  there is a $\gamma$ in $F_n$ such that
$\Delta(r_\gamma ,r^i)\geq |C^i \cap M_n|$,
$\gamma$ is in the domain of $q_n$, and $q_n(\gamma) = \sigma_i$,
\item  \label{in_Dn} $p_n \in M_n[G] \cap D_n$.
\end{enumerate}

It is not difficult to construct these sequences.
It follows directly from Clause \ref{in_Dn} that if we set 
$$G'=\{s \in N_0[G]: (\exists n \in \omega) q_n \leq s\}$$
then $\dot q[G] \in G'$ and $G'$ is $\dot \Qcal$-generic over $N_0[G]$.
Notice that, for every $i,n \in\omega$, $q_n$ is consistent with
$(C^i,r^i)\mapsto\sigma_i$.
Define a name $\bar s$ by 
$$t\Vdash \text{``}\bar s(\beta) =\begin{cases}

q_n(\beta) & \text{ if } \beta\in\dom(q_n),\\
\sigma_i& \text{ if } \beta=\delta \text{ and } t\Vdash 
\text{``} \dot C_\delta =C^i\text{ and } \dot r_\delta=r^{i}\text{''}.
\end{cases}$$
It is easy to see that if $t$ is a lower bound for $G$ and is
$\Pcal$-generic over $N_1$ and $N_2$ then
$t\Vdash `` \dot s  \in\dot \Qcal$''
and obviously $\dot s$ will be lower bound for $G'$.

\end{proof}

\begin{proof} (of Theorem \ref{decomp_Pcr}) Let $\Vbf$ be a model of $\CH$.
Construct a countable support iteration
$\Pcal_{\omega_2} =\<\Pcal_\alpha,\dot\Qcal_\alpha:\alpha<\omega_2\>$ 
such that for every $\alpha<\omega_2$ we have
$\Vdash_\alpha$ ``$\dot\Qcal_\alpha =\Qcal_{\vec{C},\vec{r}}$ for
some pair $\vec{C},\vec{r}$''.
Since $\CH$ holds in $\Vbf$ and $\Vdash_\alpha$ ``$|\Qcal_\alpha|=\aleph_1$''
it follows that $\Pcal_{\omega_2}$ satisfies the $\omega_2$-c.c..
A standard bookkeeping argument ensures that in $\Vbf^{\Pcal_{\omega_2}}$
every pair $\vec{C},\vec{r}$ 
admits a decomposition of $\omega_1=\bigcup_{n\in\omega} X_n$ such
that  
$$|C_\gamma\cap C_\delta| \leq \Delta(r_\gamma,r_\delta)$$ whenever
$\gamma < \delta$ are in the same $X_n$. 
By Theorem \ref{iter_tod} and Lemmas \ref{<w_1-proper_Pcr} and
\ref{2-complete_Pcr}, $\CH$ also holds in $\Vbf^{\Pcal_{\omega_2}}$
so the proof of Theorem \ref{decomp_Pcr} is complete.

\end{proof}

We will now finish the proof of Theorem \ref{guessing}.
Start with the model of the Theorem \ref{decomp_Pcr}.
First we will need the following lemma.

\begin{lem}\label{neat_ladder}
There is a ladder system $C_\delta$ indexed by the positive
countable limit ordinals
such that $C_\gamma \cap C_\delta$ is an initial segment of both
$C_\delta$ and $C_\gamma$ whenever $\gamma < \delta$ are limits.
\end{lem}

\begin{proof}
Let $h:\omega^{<\omega} \leftrightarrow \omega$ be a bijection
which satisfies $h(s) < h(t)$ whenever $s$ is an initial part of $t$.
For a fixed limit $\delta > 0$, we shall build an increasing $\omega$-sequence
$\bar \delta_n$ $(n \in \omega)$ cofinal in
$\delta$ such that for every $n$ the ordinal $\bar \delta_n$ has the form
$$\bar \delta_n = \xi + h(\<e_{\xi+\omega}^{-1}(\bar \delta_i):i<n\>)$$
for some limit ordinal $\xi$ (note that $\xi$ depends on $n$,
is possibly equal to 0 and
that this decomposition of $\delta_n$
is unique for any given $n$).

To see that this can be done, first note that
if $\delta = \xi + \omega$ for some limit ordinal $\xi$ then
$$\bar \delta_n = \xi + h(\<e_{\xi+\omega}^{-1}(\bar \delta_i):i<n\>)$$
recursively defines the sequence of $\bar \delta_n$'s.
If $\delta$ is a limit of limits, then first choose an increasing
sequence of limits
$\xi_n$ $(n < \omega)$ which is cofinal in $\delta$.
Again
$$\bar \delta_n = \xi_n + h(\<e_{\xi_n+\omega}^{-1}(\bar \delta_i):i<n\>)$$
recursively defines the sequence of $\bar \delta_n$'s.

Now suppose that for some positive limit ordinals
$\delta,\epsilon < \omega_1$ and some $m,n < \omega$
$\bar \delta_m = \bar \epsilon_n$.
We need to show that $m = n$ and that if $i < m$ then
$\bar \delta_i = \bar \epsilon_i$.
Find a unique limit ordinal $\xi$ and a unique element $t$ in $\omega^{<\omega}$
such that
$$\bar \delta_m = \xi + h(t) = \bar \epsilon_n.$$
Now notice that $m = |t| = n$ and
$$\bar \delta_i = e_{\bar \delta_m}(t(i)) = 
e_{\bar \epsilon_n}(t(i)) = \bar \epsilon_i$$
for any $i < n$.
\end{proof}

Fix $C_\delta = \{\bar \delta_n: n \in \omega\}$ as in Lemma \ref{neat_ladder}.
For simplicity, identify $\Rbb$ with $2^{\omega}$.
The domain of $F$ will consist of a countable sequence
$\vec{t} = \<t_n : n \in \omega\>$ of
functions from $\alpha$ to $2$ for some $\alpha \in \omega_1$.
Let the $n\Th$ element of the sequence
$F(\vec{t})$ be given by
$$k \mapsto t_n(\bar \delta_k)$$
where $\delta = |t|$.
Now suppose that $g:\omega_1 \to \Rbb^\omega$ is given.
For each $i$, let $\lim(\omega_1) = \bigcup_{j=0}^\infty X_{i,j}$
such that for all $\gamma < \delta$ in $X_{i,j}$
$$|C_\gamma\cap C_\delta| \leq \Delta(g(\gamma)(i),g(\delta)(i)).$$
Now it is possible to choose $f_n:\omega_1 \to 2$
in such a way that if $\delta$ is in $X_{i,j}$ then
$g(\delta)(i)$ is the mapping
$$k \mapsto f_{2^i3^j}(\bar \delta_k).$$
Thus for all limit $\delta$
the range of $g(\delta)$ is contained in the range of
$F(\vec{f} \restriction \delta)$.

\begin{remark}
Shelah has shown that $\diamond(3,=)$ is not a consequence of
$\CH$ (Section VIII.4 of \cite{proper_forcing}) and 
Eisworth  has shown
that $\diamond([\omega]^2,\omega,\not \ni)$ is not
a consequence of $\CH$ \cite{anti_dowker_ladders}.
\end{remark}

\end{document}